\documentclass[final,1p,times]{elsarticle}

\usepackage{amssymb}
\usepackage{amsthm}
\usepackage{amscd}
\usepackage{amsmath}
\usepackage{amsfonts}
\usepackage{amssymb}
\usepackage{graphicx}
\usepackage{color}
\usepackage{chngcntr}
\usepackage[pdftex]{pict2e}
\newtheorem{theorem}{Theorem}[section]
\newtheorem{proposition}[theorem]{Proposition}
\newtheorem{lemma}[theorem]{Lemma}

\newtheorem{remark}[theorem]{Remark}
\usepackage{mathrsfs}
\usepackage{titletoc}
\biboptions{sort&compress}
\usepackage{geometry}

\newcommand{\D}{\Delta}

\newcommand{\p}{\partial}

\newcommand{\be}{\begin{equation}}
	
	\newcommand{\ee}{\end{equation}}
\newcommand{\bea}{\begin{eqnarray}}
	\newcommand{\eea}{\end{eqnarray}}
\newcommand{\bna}{\begin{eqnarray*}}
	\newcommand{\ena}{\end{eqnarray*}}

\renewcommand{\O}{\Omega}
\renewcommand{\le}{\left}
\newcommand{\ri}{\right}

\makeatletter

\newcommand{\Rmnum}[1]{\expandafter\@slowromancap\romannumeral #1@}
\makeatother

\usepackage{bbding}

\journal{***}

\begin{document}
	\setlength{\unitlength}{0.20mm}
	\begin{frontmatter}

		\title{Fractional mean field equations on finite graphs}
		
		\author{Yang Liu}
		\ead{dliuyang@ruc.edu.cn}
		
		\address{School of Mathematics, Renmin University of China, Beijing, 100872, China}

		\begin{abstract}
			In this paper, the author considers the fractional mean field equation on a finite graph $G=(V,E)$, say
			\begin{equation*}
				(-\Delta)^s u=\rho\left(\dfrac{he^u}{\int_V he^ud\mu}-\dfrac{1}{|V|}\right),\quad\forall\,x\in V,
			\end{equation*}
			where $s\in(0,\,1)$, $\rho\in(-\infty,\,0)\cup(0,\,+\infty)$ are some fixed parameters, $h$ denotes a given real value function on $V$. Based on the sign of the prescribed function $h$, via the variational method, topological degree and two mean field type heat flows, the author obtains the existence of solutions for the above problem in three cases respectively. These results extend the relevant research of Lin-Yang (Calc. Var., 2021), Sun-Wang (Adv. Math., 2022) and Liu-Zhang (J. Math. Anal. Appl., 2023) in the case of $s=1$.
		\end{abstract}
		
		
	\end{frontmatter}
	
	\section{Introduction}
	
	Let $(\Sigma,g)$ be a closed Riemann surface without boundary and ${\rm dim}\Sigma\geq 1$, $\widetilde{g}=e^{2u}g$,  $\kappa$ and $K$ be the Gaussian curvatures associated to the Riemannian metrics ${g}$ and $\widetilde{g}$. Then $u$ satisfies
	\begin{equation}\label{eq:kw0}
		\Delta_g u+\kappa-K e^{2u}=0,
	\end{equation}
	where $\Delta_g$ denotes the Laplace-Beltrami operator of $g$.
	Let $v$ be a solution of $\Delta_gv=\overline{\kappa}-\kappa$, and $f=2(u-v)$, where $\overline{\kappa}$ is the integral average of
	$\kappa$ on $\Sigma$. Then (\ref{eq:kw0}) turns out to be a more common form
	$$\Delta_g f+2\overline{\kappa}-(2Ke^{2v})e^{f}=0.$$
	As it is conceivable, one can free \eqref{eq:kw0} from the geometric situation, and study
	\begin{equation}\label{eq:1}
		\Delta_g f+c-h e^{f}=0.
	\end{equation}
	This is the Kazdan-Warner equation, which comes from the basic geometric problem of prescribed Gaussian curvature \cite{C-D,C-Y,K-W-1}, and as well as appears in various contexts such as the abelian Chern-Simons-Higgs models in Physics \cite{D-J-L-W-2,N,R-T}. 
	
	More in general, in the presence of $c\not=0$, the Kazdan-Warner equation \eqref{eq:1} can be reduced to the following mean field equation
	\begin{equation}\label{eq:2}
		-\Delta_g u=\rho\le(\frac{he^u}{\int_\Sigma he^udv_g}-\frac{1}{\int_\Sigma 1dv_g}\ri),
	\end{equation}
	where $\rho>0$. In the past several decades, such equation has received high attention from mathematicians. For $\rho\in(0,8\pi)$, under some given geometric conditions, Ding-Jost-Li-Wang \cite{D-J-L-W-1} proved that \eqref{eq:2} has a solution if $0<h\in C^\infty(\Sigma)$. The existence of solutions was also obtained by Nolasco-Tarantello \cite{N-T}, and was generalized by Yang-Zhu \cite{Y-Z} to the case $0\leq h\not\equiv0$. For $\rho\in(8\pi,4\pi^2)$, under the assumptions $h\equiv1$ and $\Sigma$ is a flat torus, utilizing a min-max scheme, Struwe-Tarantello \cite{S-T} proved that \eqref{eq:2} admits a non-constant solution. For $\rho\in(8\pi, 16\pi)$, by a similar method, Ding-Jost-Li-Wang \cite{D-J-L-W-3} obtained a non-minimal solution of the equation \eqref{eq:2} when $h>0$ and the genus of $\Sigma$ is greater than zero. For $\rho\in(8k\pi,8(k+1)\pi)$ $(k\in\mathbb{N})$, Chen-Lin \cite{C-L} obtained that the Leray-Schauder degree is
	\begin{equation*}
		d_\rho=\begin{pmatrix}
			k-\chi(\Sigma)\\
			k
		\end{pmatrix},
	\end{equation*} 
	where $\chi(\Sigma)$ is the Euler number of $\Sigma$. In particular, if the genus of $\Sigma$ is positive and $h>0$, then the mean field equation \eqref{eq:2} always possesses at least one solution. After that, this result was extended by Djadli \cite{D}, Malchiodi \cite{M} to the case $\Sigma$ is a general Riemann surface. 
	
	Recently, Cast\'{e}ras \cite{C1,C2} introduced a mean field type heat flow, namely
	\begin{equation}\label{108}\left\{\begin{array}{lll}
			\dfrac{\partial}{\partial t}e^u=\Delta_g u+\rho\dfrac{he^u}{\int_\Sigma he^udv_g}-Q,\quad\left(x,\, t\right)\in \Sigma\times\left(0,\, +\infty\right),\\[3ex]
			u\left(x,\, 0\right)=u_0\left(x\right),\quad\quad\quad\quad\quad\quad\quad\quad x\in \Sigma,\end{array}\right.
	\end{equation}
	where $\rho>0$, $0<h\in C^\infty(\Sigma)$, $u_0\in C^{2+\alpha}(\Sigma)$ is the initial function with $\alpha\in(0,1)$, $Q\in C^{\infty}(\Sigma)$ is a given function satisfying $\int_\Sigma Qdv_g=\rho$. For $\rho\not=8m\pi$ $(m\in\mathbb{Z}^+)$, provided with an appropriate initial function, Cast\'{e}ras proved that there exists a unique solution $u\in C^{2+\alpha,1+\frac{\alpha}{2}}_{loc}(\Sigma\times[0,+\infty))$ of the flow \eqref{108}; even more, he showed that $u$ converges in $H^2(\Sigma)$ to a function $u_\infty\in C^\infty(\Sigma)$, the solution of the following mean field equation
	\begin{equation*}
		-\Delta_gu_\infty=\rho\dfrac{he^{u_\infty}}{\int_\Sigma he^{u_\infty} dv_g}-Q.
	\end{equation*}
	In this spirit, it was supplemented by Li-Zhu \cite{LZ} that the flow also converges at a critical case $\rho=8\pi$. The same conclusion was as well as generalized by Sun-Zhu \cite{S-Z1} to the case $0\leq h\not\equiv0$. By constructing a non-local gradient-like flow, the existence of solutions for the equation \eqref{eq:2} was also drawn by Li-Xu \cite{L-X} when $h$ changes sign and $\rho=8\pi$. Very recently, using an elliptic method, Wang-Yang \cite{W-Y} derived the existence of solutions for the \textbf{G}-invariant flow \eqref{108} in the case $h$ changes sign, where $\textbf{G}=\{\sigma_1,\,\cdots,\,\sigma_N\}$  is a finite isometric group acting on a closed Riemann surface. For recent progress in this direction, we refer readers to \cite{W-Y1,W-Y2,Y-Y,Z} for more details.\\
	
	In recent years, the research of discrete partial differential equations on graphs, generated in image processing and artificial intelligence applications, has become a highly focused direction in the literature. Since then many efforts have been made to study the existence of the mean field equation \eqref{eq:2} in discrete setting. Utilizing the variational method, Zhu \cite{Z2} studied the mean field equations for the equilibrium turbulence and Toda systems in the case of $0\leq h\not\equiv0$ and $\rho>0$ on a connected finite graph. Then, when $\rho\in(-\infty,0)\cup(0,+\infty)$ and $h>0$, Liu \cite{L} considered the discrete form of \eqref{eq:2} on a finite graph, say
	\begin{equation}\label{03}
		-\Delta u=\rho\left(\dfrac{he^u}{\int_V he^ud\mu}-\dfrac{1}{|V|}\right),\quad\forall\,x\in V,
	\end{equation}
	where $-\D$ is the Laplacian on graphs. By calculating the Brouwer degree of relevant map, Liu pointed out that the equation \eqref{03}  admits at least one solution. More recently, the heat flow method of Cast\'{e}ras \cite{C1,C2} was extended by Lin-Yang \cite{ly0} from a closed Riemann surface to a finite graph. They investigated the mean field type flow \eqref{108} on finite graphs in the case of $h\equiv1$, namely
	\begin{equation}\label{04}\left\{\begin{array}{lll}
			\dfrac{\p}{\p t}\phi\left(u\right)=\Delta u-Q+\rho\dfrac{e^u}{\int_V e^ud\mu},\quad\left(x,\, t\right)\,\in \, V\times\left(0,\, +\infty\right),\\[3ex]
			u\left(x,\, 0\right)=u_0\left(x\right),\quad\quad\quad\quad\quad\quad\ \ x\,\in \,V,\end{array}\right.
	\end{equation}
	where $\rho\in\mathbb{R}$, $Q$ is a function with $\int_V Qd\mu=\rho$ and $\phi:\mathbb{R}\rightarrow\mathbb{R}$ is a $C^1$ function satisfying
	$$\lim_{s\rightarrow-\infty}\phi\left(s\right)=0,\quad \phi^\prime\left(s\right)>0\ \ \mathrm{for\ all}\ \ s\in\mathbb{R},\quad\ \inf_{s\in[0,+\infty )}\phi^\prime\left(s\right)>0.$$
	In particular, they derived that for any $\rho\in\mathbb{R}$, the solution $u(t)$ of flow \eqref{04} converges to a function $u_\infty$, the solution of the following mean field equation 
	\begin{equation*}
		-\Delta u_\infty+Q=\rho\dfrac{e^{u_\infty}}{\int_V e^{u_\infty}d\mu}.
	\end{equation*}
	Later on, Liu-Zhang \cite{L-Z} eliminated the restriction on $h$ and extended the flow method of Wang-Yang \cite{W-Y} to the case $h$ changes sign on finite graphs, and studied the flow 
	\begin{equation*}\left\{\begin{array}{lll}
			\dfrac{\p}{\p t}e^u=\Delta u+\rho\left(\dfrac{he^u}{\int_V he^ud\mu}-\dfrac{1}{|V|}\right),\quad\left(x,\, t\right)\,\in \, V\times\left(0,\, +\infty\right),\\[3ex]
			u\left(x,\, 0\right)=u_0\left(x\right),\quad\quad\quad\quad\quad\quad\quad\ \  x\,\in \,V,\end{array}\right.
	\end{equation*}
	where $\rho$ is a positive real number, supplementing the results of Lin-Yang \cite{ly0}. For other discrete partial differential equations on graphs, it has been studied extensively so far, see papers, for example, the Schr\"{o}dinger equations \cite{A-Y-Y-3,Z-Z}, the Kazdan-Warner equations \cite{C-M,G,A-Y-Y-2,K-S,L-S}, the Chern-Simons-Higgs equations \cite{HS,H-L-Y,H-W-Y} and so on  \cite{B-S-W,H-S-Z,L-W,P-S,W}, and the references therein.
	
	Regarding the latest progress related to this topic, it is worth noting that the analysis on some more specific discrete spaces, such as Cayley graphs and lattice graphs, has been posed by Hua-Li-Wang-Xu \cite{H-L,H-L-W,H-X}. Moreover, Gu-Huang-Sun-Lu \cite{G-H-S,LS} investigated the problem of semilinear elliptic inequality on graphs. Recently, Shao-Yang-Zhao \cite{S-Y-Z} have developed the Sobolev spaces theory on locally finite graphs and obtained the reflexivity, completeness, separability, and Sobolev inequalities of $W^{m,p}(V)$. More recently, we would like to also mention that Shao-Tian-Zhao \cite{S-T-Z} established a coherent framework of variational methods on hypergraphs, including the propositions of calculus and the theory of function spaces on hypergraphs. As for other related works on graphs, there are more recent progresses, see for instance \cite{C-W-Y,C-H,I-B-R,L-Z2,S}.\\
	
	Let us come back to introduce the fractional Laplace operator in Euclidean space. As an important example of a nonlocal pseudo differential operator, the fractional Laplacian and its corresponding nonlinear problems have impressive applications in optimization, conservation laws, stratified materials, gradient potential theory and so on, see papers \cite{B-K,Mg,S-V,D-L}. We recommend several books or papers \cite{A,R-S,D-P-V,K,C-K-S,S-V1} about classic fractional Laplace operator theory to readers, among which Kwa\'{s}nicki \cite{K} introduced ten equivalent definitions of the fractional Laplace operator on a continuous function space. Specifically, one of the definitions is given by the heat kernel, namely, for any fixed $s\in(0,\,1)$, and $u\in C^\infty(\mathbb{R}^N)$
	\begin{equation}\label{10}
		(-\D)^su(x)=\frac{s}{\Gamma(1-s)}\int_0^{+\infty}\le(u(x)-e^{t\D}u(x)\ri)t^{-1-s}dt\in L^2(\mathbb{R}^N).
	\end{equation}
	Here $\Gamma(\cdot)$ represents the Gamma function, $e^{t\D}$ denotes the heat semigroup of the Laplace operator $-\D$. Nowadays, a great interest has emerged in the discrete fractional Laplace operator and related fractional problems. The definition of $(-\D)^s$ on a lattice graph $\mathbb{Z}^d$ is given by Lizama-Murillo-Arcila \cite{L-M} and Wang \cite{Wl} through discrete Fourier transform, by Ciaurri et al. \cite{C-G-R-T-V,C-R-S-T-V}, Lizama-Roncal \cite{L-R} and Stinga-Torrea \cite{S-T1} with the semigroup method. Recently, Wang \cite{Wj} and Zhang-Lin-Yang \cite{Z-L-Y} provided definitions of $(-\D)^s$ on locally finite graphs and finite graphs by the heat semigroup respectively. In addition, for any positive real number $s$, Zhang-Lin-Yang \cite{Z-L-Y} studied the fractional Kazdan-Warner equation on a finite graph
	\begin{equation}\label{05}
		(-\D)^su=\kappa e^u-c,\quad\text{in}\ V.
	\end{equation}
	Via variational principles and the method of upper and lower solutions, they obtained several existence results. \\
	
	In this paper, our aim is to study the mean field equation \eqref{03} involving the discrete fractional Laplace operator $(-\D)^s$ on a finite graph. Except for the fractional Kazdan-Warner equation \eqref{05} treated by Zhang-Lin-Yang \cite{Z-L-Y}, little is known regarding the fractional mean field equation. It is impossible to employ the approaches in \cite{Z-L-Y} to investigate the existence of the solutions. According to the sign of $h$, we will utilize different methods, such as the variational method, topological degree and two mean field type heat flows, to discuss the existence of solutions separately. Compared with all the existing relevant results, the innovation of our present works mainly lies in
	
	\noindent(a) Fractional Laplacian is a class of non-local operator. This is the essential difference from the counterparts in \cite{Z2,L,ly0,L-Z}, where the Laplace operator was considered.\\
	(b) For any $s\in(0,\,1)$, we completely solve the problem of the existence of solutions, regardless of the sign of $h$.\\
	(c) In the graph setting, $\rho$ can take any value in $(-\infty,0)\cup(0,+\infty)$. However, it is clear to see that $\rho$ is assumed to be greater than $0$ on compact Riemann surfaces.
	
	The remaining parts of this paper are organized as follows: In Section 2, we introduce the discrete version of the fractional Laplacian on graphs and state our main results. In Sections 3-5, we shall look for the solutions in three cases, respectively.

	\section{Preliminaries and Main result}
	
	To start with, let us present some basic definitions on graphs. A graph $G=\left(V,E\right)$ is said to be a finite graph, if $\sharp V$ is a finite number, or $V$ only contains finitely many vertices, where ${V}$ is the vertex set and $E$ denotes the edge set. Let ${\mu:V\rightarrow\mathbb{R}^+}$ be a positive measure and $w:E\rightarrow \mathbb{R}^+$ be a positive symmetric weight satisfying ${w_{xy}>0}$ and ${w_{xy}=w_{yx}}$ for every edge ${xy\in E}$. Throughout this paper, we always assume that $G=(V,E)$ is a finite graph, which satisfies\\
	
	\noindent(a) (Simple) $G$ contains neither loops nor multiple edges.\\
	(b) (Undirected) The edges of $G$ are directionless, which means $x$ and $y$ can reach each other for any $xy\in E$.\\
	(c) (Connected) For any $x,y\in V$, there exist finite edges connecting $x$ and $y$.\\
	
	Let us define $C(V):=\{f:V\rightarrow\mathbb{R}\}$. The integral of a function $u\in C(V)$ is represented by
	\begin{equation*}
		\int_{V}ud\mu=\sum_{x\in V}\mu(x)u(x).
	\end{equation*}
	For any $1\leq p\leq+\infty$, we consider a function space 
	$${\ell^p(V)}=\left\{u\in C(V):\|u\|_{\ell^p(\O)}<+\infty\right\},$$
	endowed with a norm
	\begin{equation}\label{22}
		\|u\|_{\ell^p(V)}= \left\{\begin{aligned}
			&\left(\sum_{x\in V}\mu(x)|u(x)|^p\right)^{\frac{1}{p}}, &1\leq p<+\infty,\\
			&\max_{x\in V}|u(x)|, &p=+\infty.\ \ \ \ \ \ \,\end{aligned}\right.
	\end{equation}
	For convenience, we will write $\|\cdot\|_p$ and $\|\cdot\|_\infty$ instead of $\|\cdot\|_{\ell^p(V)}$ and $\|\cdot\|_{\ell^\infty(V)}$, when there is no confusion.
	
	\subsection{The fractional Laplacian on graphs}
	
	Before introduing the fractional Laplace operator on graphs formally, let us recall some known concepts of Laplace operator on graphs. The Laplacian on $G$ is defined as
	\begin{equation}\label{9}
		\Delta u(x)=\frac{1}{\mu(x)}\sum_{y\sim x}w_{xy}(u(y)-u(x)),
	\end{equation}
	where $y\sim x$ means $y$ is adjacent to $x$, i.e, $xy\in E$. Clearly, $\D$ is a local operator. Without loss of generality, we may assume that $V=\{x_1,\,\cdots,\, x_n\}$ and $\sharp V=n$ for some positive integer $n$. Notice that $-\D$ is a nonnegative definite symmetric operator. Its eigenvalues are written as $0=\lambda_0<\lambda_1\leq\lambda_2\leq\cdots\leq\lambda_{n-1}$,
	and the corresponding eigenfunction spaces are one dimension. Let the corresponding eigenvectors be $\phi_0,\,\phi_1,\,\ldots,\,\phi_{n-1}$, and set these to be orthonormal. That is, for any $i,j=0,\,1,\,\ldots,\,n-1$,
	\begin{equation*}
		\langle\phi_i,\phi_j\rangle=\sum_{x\in V}\mu(x)\phi_i(x)\phi_j(x)=\left\{\begin{aligned}
			1,\quad i=j,\\
			0,\quad i\not=j.\end{aligned}\right.
	\end{equation*}
	 Here the orthonormality is the inner product associated to the measure $\mu$. 

	In what follows, the fractional Laplace operator $(-\D)^s$ can be explicitly expressed on a finite graph according to the previous works by Zhang-Lin-Yang \cite{Z-L-Y}, which is the graph case of Bochner’s definition in the Euclidean setting by Kwa\'{s}nicki \cite{K}. Let $s\in(0,\,1)$, the fractional Laplace operator $(-\D)^s$ is defined as
		\begin{equation}\label{7}
			(-\D)^su(x)=\frac{1}{\mu(x)}\sum_{y\in V, y\not=x}W_s(x,y)(u(x)-u(y)),
		\end{equation}
		with a function
		\begin{equation}\label{8}
			W_s(x,y)=-\mu(x)\mu(y)\sum_{i=0}^{n-1}\lambda_i^s\phi_i(x)\phi_i(y),\quad\forall\ x\not=y\in V,
		\end{equation}
	where $n=\sharp V\geq2$, $\lambda_i$ is the eigenvalue of the Laplacian $-\D$, and $\phi_i$ is the corresponding orthonormal eigenfunction, the function $W_s(x,y)$ satisfies
	$$0<W_s(x,y)=W_s(y,x)<+\infty,\quad\forall\ x\not=y\in V.$$
	 It is obvious that the fractional Laplacian is well-defined for any $u\in C(V)$.
	
	Now we define the fractional gradient form of $u\in C(V)$ at $x\in V$. With no loss of generality, we may assume $V=\{x,\,y_1,\,y_2,\,\ldots,\,y_{n-1}\}$ with $\sharp V=n\geq2$. Then the fractional gradient of $u$ reads as
	\begin{equation*}
		\nabla^su(x)=\le(\sqrt{\frac{W_s(x,y_1)}{2\mu(x)}}(u(x)-u(y_1)),\,\cdots,\,\sqrt{\frac{W_s(x,y_{n-1})}{2\mu(x)}}(u(x)-u(y_{n-1}))\ri)\in\mathbb{R}^{n-1},
	\end{equation*}
	where $W_s(x,y_i)$ is defined by \eqref{8}. The inner product of two gradients $\nabla^su$ and $\nabla^sv$ at $x\in V$ is 
	\begin{equation}\label{13}
		\nabla^s u\nabla^s v(x)=\frac{1}{2\mu(x)}\sum_{y\in V, y\not=x}W_s(x,y)(u(x)-u(y))(v(x)-v(y)),
	\end{equation}
	and the norm of $\nabla^s u$ is given as
	\begin{equation}\label{19}
		|\nabla^s u|(x)=\sqrt{\nabla^s u\nabla^s u(x)}=\le(\frac{1}{2\mu(x)}\sum_{y\in V, y\not=x}W_s(x,y)(u(x)-u(y))^2\ri)^{\frac{1}{2}}.
	\end{equation}
	For any $s\in(0,\,1)$, the fractional Sobolev space $W^{s,2}(V)$ is defined by
	\begin{equation*}
		W^{s,2}(V)=\le\{u\in C(V):\int_{V}\le(|\nabla^su|^2+|u|^2\ri)d\mu<+\infty\ri\},
	\end{equation*}
	with a norm
	\begin{equation}\label{73}
		\|u\|_{W^{s,2}(V)}=\le(\int_{V}\le(|\nabla^su|^2+|u|^2\ri)d\mu\ri)^{\frac{1}{2}}.
	\end{equation}
	It is clear that $W^{s,2}(V)$ is a Hilbert space with the inner product
	\begin{equation*}
		\langle u,\,v\rangle_{W^{s,2}(V)}=\int_{V}\le(\nabla^su\nabla^sv+uv\ri)d\mu.
	\end{equation*}
	Since $V$ only contains finitely many vertices, it is possible to see that $W^{s,2}(V)=C(V)$ is exactly $\mathbb{R}^{n}$, a finite dimensional linear space. 

	Despite some similarities, by means of \eqref{7} and \eqref{8}, the nonlocal linear operator $(-\D)^s$ is very different from $-\D$ given in \eqref{9}; even more, it is interesting that the fractional Laplacian $(-\D)^s$ on a finite graph can be fully represented by the eigenvalues of the Laplacian $-\D$ and the corresponding orthonormal eigenfunction. However, $(-\D)^s$ on infinite graph does not have such a representation. It is expected that, corresponding to the eigenvalue $\lambda^s_i$ of the fractional Laplacian $(-\D)^s$, $\phi_i$ remains an orthonormal eigenfunction (see \cite{Z-L-Y}, Proposition 3.3). Namely, for any $x\in V$, there holds
		\begin{equation}\label{06}
			(-\D)^s\phi_{i}(x)=\lambda_i^s\phi_{i}(x),\quad\forall\ i=0,\,1,\,\ldots,\,n-1,
		\end{equation}
	where $0=\lambda_0<\lambda_1\leq\lambda_2\leq\cdots\leq\lambda_{n-1}$ are eigenvalues of $-\D$, and $\phi_0,\,\phi_1,\,\ldots,\,\phi_{n-1}$ are the corresponding orthonormal eigenfunctions.
	
	Finally in this subsection, we highlight the formula of integration by parts on graphs, which plays a crucial role in the analysis on graphs. Here we provide a direct proof different from Zhang-Lin-Yang (\cite{Z-L-Y}, Proposition 3.4).
	\begin{lemma}\label{l2}
		(Formula of integration by parts) Let $s\in(0,\,1)$. Suppose that $u\in W^{s,2}(V)$ and $\nabla^s u\nabla^s v$ is well defined as \eqref{13}. Then we have
		\begin{equation}\label{12}
			\int_Vu(-\D)^s v d\mu=\int_{V}\nabla^s u\nabla^s vd\mu=\int_Vv(-\D)^s u d\mu,\quad\forall\ v\in W^{s,2}\left(V\right),
		\end{equation}
		where $(-\D)^s$ is given by \eqref{7}.
	\end{lemma}
	\begin{proof}
		Note that $W_s(x,y)$ is a symmetric function for any $x\not=y\in V$. It follows from the definition of associated gradient in \eqref{13} and a straightforward calculation that
		\begin{align}\label{14}
			\nonumber\int_{V}\nabla^s u\nabla^s vd\mu&=\frac{1}{2}\sum_{x\in V}\sum_{y\in V, y\not=x}W_s(x,y)(u(x)-u(y))(v(x)-v(y))\\
			\nonumber&=\frac{1}{2}\sum_{x\in V}\sum_{y\in V, y\not=x}W_s(x,y)(u(x)-u(y))v(x)-\frac{1}{2}\sum_{x\in V}\sum_{y\in V, y\not=x}W_s(x,y)(u(x)-u(y))v(y)\\
			\nonumber&=\frac{1}{2}\sum_{x\in V}\sum_{y\in V, y\not=x}W_s(x,y)(u(x)-u(y))v(x)-\frac{1}{2}\sum_{y\in V}\sum_{x\in V, x\not=y}W_s(x,y)(u(x)-u(y))v(y)\\
			\nonumber&=\sum_{x\in V}\sum_{y\in V, y\not=x}W_s(x,y)(u(x)-u(y))v(x)\\
			&=\int_Vv(-\D)^s u d\mu.
		\end{align}
		Similarly as we did in \eqref{14}, we can also obtain that
		$$\int_{V}\nabla^s u\nabla^s vd\mu=\int_Vu(-\D)^s v d\mu,$$
		which gives the desired result.
	\end{proof}
	
	\subsection{Main results}
	
	In the present paper, we mainly consider the existence of solutions for the fractional mean field equation on a finite graph $G=(V,E)$, say
	\begin{equation}\label{15}
		(-\Delta)^s u=\rho\left(\dfrac{he^u}{\int_V he^ud\mu}-\dfrac{1}{|V|}\right),\quad\forall\,x\in V,
	\end{equation}
	where $s\in(0,\,1)$, $\rho\in(-\infty,\,0)\cup(0,\,+\infty)$ are some fixed parameters, $h\in C(V)$ stands for a given real value function on $V$, $|V|=\sum_{x\in V}\mu(x)$ means the volume of $V$. Let us define two closed nonempty subspaces of $W^{s,2}(V)$ by
	\begin{equation}\label{16}
		\mathscr{H}=\le\{u\in W^{s,2}(V):\int_{V}ud\mu=0\ri\}
	\end{equation}
	and
	\begin{equation}\label{35}
		\mathscr{M}=\le\{u\in W^{s,2}(V):\int_{V}he^ud\mu=1\ri\}.
	\end{equation}
	
	If $h$ is a positive function on $V$, we concern the variational functional related to \eqref{15}. Define $J_\rho:\mathscr{H}\rightarrow\mathbb{R}$ by
	\begin{equation}\label{17}
		J_\rho(u)=\frac{1}{2}\int_{V}|\nabla^su|^2d\mu-\rho\log\int_{V}he^ud\mu.
	\end{equation}
	By employing the direct variational method, we reach the first result.
	\begin{theorem}\label{t1}
		Let $G=(V,E)$ be a connected finite graph, $\mathscr{H}$ and $J_\rho$ be defined as \eqref{16} and \eqref{17} respectively. Suppose that $h\in C(V)$ is a positive function, then for any $\rho\in(-\infty,\,0)\cup(0,\,+\infty)$ and $s\in(0,\,1)$,  $J_\rho$ has a minimizer $u\in\mathscr{H}$. Moreover, the mean field equation \eqref{15} has a solution $u\in\mathscr{H}$.
	\end{theorem}
	
	If $h\in C(V)$ is a nonnegative function satisfying $0\leq h(x)\not\equiv0$ for all $x\in V$, using the method of topological degree, we shall obtain the existence results of $\eqref{15}$, which is our second result.
	\begin{theorem}\label{t2}
		Let $G=(V,E)$ be a connected finite graph, $\mathscr{M}$ be given by \eqref{35}. Suppose that $h\in C(V)$ satisfying $0\leq h\not\equiv0$, then for any $\rho\in(-\infty,\,0)\cup(0,\,+\infty)$ and $s\in(0,\,1)$, the mean field equation \eqref{15} has at least one solution $u\in\mathscr{M}$.
	\end{theorem}
	
	The topological degree in the graph setting was first pioneered by Sun-Wang \cite{S-w} to the Kazdan-Warner equation, then adapted to  the mean field equation by Liu \cite{L}, and to the Chern-Simons-Higgs models by Li-Sun-Yang \cite{L-S-Y}. In the proof of Theorem \ref{t2}, the first step of this approach is to establish a priori estimate. This closely relies on the structure of the set $\{x\in V: h(x)=0\}$, which is essentially different from the counterparts in \cite{L,S-w,L-S-Y}. Once the uniform estimate is obtained, the topological degree follows and can be employed to find the solutions of \eqref{15}. This is done using the homotopy invariance and Kronecker's existence Theorem. 
	
	If $h\in C(V)$ is a sign-changing prescribed function satisfying $\max_{x\in V}h(x)>0$, the problem becomes more difficult. Following the lines of Wang-Yang \cite{W-Y} in the Euclidean setting, we wish to get rid of the non-negativity assumption of $h$. Here, we overcome the defect by observing the convergence results of two different heat flows. The first heat flow to be considered in this paper can be stated as follows: 
	\begin{equation*}\left\{\begin{array}{lll}
			\dfrac{\p}{\p t}e^u=-(-\Delta)^s u+\rho\left(\dfrac{he^u}{\int_V he^ud\mu}-\dfrac{1}{|V|}\right),\quad\left(x,\, t\right)\,\in \, V\times\left(0,\, +\infty\right),\\[3ex]
			u\left(x,\, 0\right)=u_0\left(x\right),\quad\quad\quad\quad\quad\quad\quad\ \  x\,\in \,V,\end{array}\right.
	\end{equation*}
	where $\rho$ is a positive real number. Motived by our recent work  \cite{L-Z}, we have the following results without the crucial assumption on the prescribed function $h$.
	\begin{theorem}\label{t3}
		Let $G=\left(V, \, E\right)$ be a connected finite graph. Suppose that $\rho>0$ is a given parameter, and that $h\in C(V)$ changes sign satisfying $\max_{V}h>0$. Then for any initial function $u_0\in C(V)$, it results:\\
		(i) there exists a unique solution $u:V\times[0,\, +\infty)\rightarrow\mathbb{R}$ of the first mean field type flow;\\
		(ii) there exists some function $u_\infty$ such that $u\left(x,\, t\right)$ converges to $u_\infty$ uniformly in $x\in V$ as $t$ tends to $+\infty$; even more, $u_\infty$ is a solution of the mean field equation (\ref{15}).
	\end{theorem}
	In proving Theorem \ref{t3}, we emphasize the crucial role played by Liu-Zhang (\cite{L-Z}, Lemma 3.3), in the study of compactness result for long time solutions. The key advantage of \cite{L-Z} is that $h$ being nonnegative is not needed. Finally, let us define the scalar curvature type function $R$ by
	\begin{equation*}
		R=\frac{1}{\rho}e^{-u}\le((-\D)^su+\frac{\rho}{|V|}\ri),
	\end{equation*}
	and the normalized coefficient $\alpha(t)$ by
	\begin{equation*}
		\alpha(t)=\frac{\int_{V}hRe^ud\mu}{\int_{V}h^2e^ud\mu}.
	\end{equation*}
	We concern the second heat flow of the form
	\begin{equation*}\left\{\begin{array}{lll}
			\dfrac{\p}{\p t}u=\alpha(t)h-R,\quad\left(x,\, t\right)\,\in \, V\times\left(0,\, +\infty\right),\\[3ex]
			u\left(x,\, 0\right)=u_0\left(x\right)\in\mathscr{M},\quad\ \  x\,\in \,V,\end{array}\right.
	\end{equation*}
	where $\rho>0$, $h\in C(V)$ changes sign and $\mathscr{M}$ is given in \eqref{35}. The second flow is a non-local gradient-like flow, which was first introduced by Li-Xu \cite{L-X} in the Euclidean setting. By studying the properties of the second flow, we can obtain the solutions of \eqref{15} in $\mathscr{M}$, which is our last result.
	
	\begin{theorem}\label{t4}
		Under the conditions of Theorem \ref{t3}. For any initial function $u_0\in\mathscr{M}$, we have:\\
		(i) there exists a unique solution $u:V\times[0,\, +\infty)\rightarrow\mathbb{R}$ of the second flow; \\
		(ii) there exists an increasing sequence $t_n\rightarrow+\infty$ and some function $u_\infty\in\mathscr{M}$ such that $u\left(x,\, t_n\right)$ converges to $u_\infty$ uniformly in $x\in V$ as $n\rightarrow+\infty$. Moreover, $u_\infty$ is a solution of the mean field equation (\ref{15}).
	\end{theorem}

	\section{The case $h>0$}
	
	Let $h$ be a given positive function on $V$. Our intention in this section is to prove Theorem \ref{t1} under the variational framework. For this purpose, we begin our study by introducing discrete version of Poincar\'{e} inequality and Moser-Trudinger inequality, which will be frequently used in the following discussions.
	
	\begin{lemma}\label{l3}
		(Poincar\'{e} inequality) Let $G=\left(V,\, E\right)$ be a connected finite graph. For any $u\in \mathscr{H}$, there exists some constant $C>0$ depending only on $G$ such that
		\begin{equation*}
			\int_V u^2d\mu\leq C(G)\int_V |\nabla^s u|^2d\mu.
		\end{equation*}
	\end{lemma}
	\begin{proof}
		Arguing by contradiction, we assume that for any $k\in\mathbb{N}$, there exists $\{u_k\}\subset W^{s,2}(V)$ so that
		\begin{equation}\label{18}
			\int_V |\nabla^s u_k|^2d\mu<\frac{1}{k},\quad\int_V u_kd\mu=0,\quad\int_V u_k^2d\mu=1.
		\end{equation}
		Hence $\{u_k\}$ is bounded in $W^{s,2}(V)$. According to the equivalent norms, it results that $\{u_k\}$ is also bounded in $\ell^\infty(V)$. Then we can find some $u_0\in W^{s,2}(V)$ such that passing to a subsequence $\{u_k\}$ (hereafter, we do not distinguish sequence and its subsequence),
		$$u_k\rightarrow u_0\quad\textrm{in}\quad W^{s,2}(V),\quad\textrm{as}\quad k\rightarrow+\infty.$$
		Taking $k\rightarrow+\infty$ in \eqref{18}, one arrives at
		\begin{equation}\label{20}
			\int_V |\nabla^s u_0|^2d\mu=0,\quad\int_V u_0d\mu=0,\quad\int_V u_0^2d\mu=1.
		\end{equation}
		Noting that \eqref{19}, we infer from the first equality in \eqref{20} that
		\begin{equation*}
			\int_V |\nabla^s u_0|^2d\mu=\frac{1}{2}\sum_{x\in V}\sum_{y\in V, y\not=x}W_s(x,y)(u_0(x)-u_0(y))^2=0,
		\end{equation*}
		which together with the second equality implies $u_0\equiv0$,  in contradiction with the third equality in \eqref{20}. Thus we complete the proof of the lemma.
	\end{proof}
	
	\begin{lemma}\label{l4}
		(Moser-Trudinger inequality) Let $G=\left(V,\, E\right)$ be a connected finite graph. For any $\beta>0$, there exists a constant $C>0$ depending only on $\beta$ and $G$ such that for all functions $v\in \mathscr{H}$ with $\int_V |\nabla^s v|^2d\mu\leq1$, there holds
		\begin{equation*}
			\int_V e^{\beta v^2}d\mu\leq C(G,\beta).
		\end{equation*}
	\end{lemma}
	\begin{proof}
		For any $v\in \mathscr{H}$ with $\int_V |\nabla^s v|^2d\mu\leq1$, by Lemma \ref{l3} one knows
		$$\int_V v^2d\mu\leq C(G)\int_V |\nabla^s v|^2d\mu\leq C(G).$$
		Then for any $x\in V$ we deduce that
		$$v^2\left(x\right)\leq\frac{1}{\min_{x\in V}\mu\left(x\right)}\int_V v^2d\mu\leq C(G).$$
		Therefore,
		$$\int_V e^{\beta v^2}d\mu=\sum_{x\in V}\mu\left(x\right)e^{\beta v^2}\leq e^{\beta C(G)}|V|\leq C(G,\beta),$$
		which gives the desired result.
	\end{proof}
	
	Since the case of $\rho<0$ is considered, next we need to estimate the lower bound of the term $\int_{V}he^ud\mu$.
	\begin{lemma}\label{l5}
		Let $G=\left(V,\, E\right)$ be a connected finite graph and $h$ be a positive function on $V$. There exist constants $C_1(G,h)>0$ and $C_2(G)<0$ such that for any $u\in\mathscr{H}$, there holds
		\begin{equation}\label{21}
			\int_{V}he^ud\mu\geq C_1e^{C_2\|\nabla^su\|_2}.
		\end{equation}
	\end{lemma}
	\begin{proof}
		Denote $\mu_0=\min_{x\in V}\mu(x)$. In view of \eqref{22}, we have
		\begin{equation*}
			\|u\|_2=\left(\sum_{x\in V}\mu(x)|u(x)|^2\right)^{\frac{1}{2}}\geq\le(\mu_0\sum_{x\in V}|u(x)|^2\ri)^{\frac{1}{2}}\geq\mu_0^{\frac{1}{2}}\max_{x\in V}|u(x)|,
		\end{equation*}
		which leads to
		\begin{equation}\label{23}
			\|u\|_\infty\leq\mu_0^{-\frac{1}{2}}\|u\|_2.
		\end{equation}
		It follows from \eqref{23} and Lemma \ref{l3} that
		\begin{align*}
			\int_{V}he^ud\mu\geq h_0|V|e^{-\|u\|_\infty}\geq h_0|V|e^{-\mu_0^{-\frac{1}{2}}\|u\|_2}\geq h_0|V|e^{-\le(\frac{C(G)}{\mu_0}\ri)^{\frac{1}{2}}\|\nabla^su\|_2},
		\end{align*}
		where $h_0=\min_{x\in V}h(x)>0$. In this case, we infer \eqref{21} with $C_1=h_0|V|>0$ and $C_2=-C(G)^\frac{1}{2}\mu_0^{-\frac{1}{2}}<0$. 
	\end{proof}
	
	Recall the variational functional  $J_\rho:\mathscr{H}\rightarrow\mathbb{R}$ associated with \eqref{15}, given by 
	\begin{equation*}
		J_\rho(u)=\frac{1}{2}\int_{V}|\nabla^su|^2d\mu-\rho\log\int_{V}he^ud\mu,
	\end{equation*}
	where
	$$\mathscr{H}=\le\{u\in W^{s,2}(V):\int_{V}ud\mu=0\ri\}.$$
	Now we are in position to prove Theorem \ref{t1} via the variational method. \\
	
	\noindent$\textbf{\emph{Proof of Theorem \ref{t1}}.}$ We divide the proof into three steps.\\
	
	\textbf{Step 1.} We first prove that the variational functional $J_\rho$ is bounded from below in $\mathscr{H}$.
	
	\noindent \textbf{Case 1:} $\rho>0$. Using the Young inequality, one can obtain for any $\epsilon>0$ and $u\not\equiv C\in\mathscr{H}$ that
	$$u\leq\epsilon\|\nabla^su\|_2^2+\frac{1}{4\epsilon}\le(\frac{u}{\|\nabla^su\|_2}\ri)^2.$$
	By Lemma \ref{l4}, we conclude 
	\begin{align}\label{24}
		\nonumber\log\int_{V}he^ud\mu&\leq\log\max_{V}h+\log\int_{V}e^{\epsilon\|\nabla^su\|_2^2+\frac{1}{4\epsilon}\le(\frac{u}{\|\nabla^su\|_2}\ri)^2}d\mu\\
		\nonumber&\leq\log\max_{V}h+\log C(G,\epsilon)+\epsilon\|\nabla^su\|_2^2\\
		&=C(G,h,\epsilon)+\epsilon\|\nabla^su\|_2^2.
	\end{align}
	Here and in the sequel, we denote by $C$ a constant which may be different from line to line. Submitting \eqref{24} into \eqref{17}, we have for any $u\in\mathscr{H}$
	\begin{align*}
		J_\rho(u)\geq\frac{1}{2}\int_{V}|\nabla^su|^2d\mu-\epsilon\rho\|\nabla^su\|_2^2-C(G,h,\epsilon,\rho).
	\end{align*}  
	Choosing $\epsilon=1/4\rho$, in particular, it implies
	\begin{equation}\label{25}
		J_\rho(u)\geq\frac{1}{4}\int_{V}|\nabla^su|^2d\mu-C(G,h,\rho).
	\end{equation}
	\noindent \textbf{Case 2:} $\rho<0$. Applying the Young inequality, reviewing \eqref{17} and \eqref{21}, we have that for any $\epsilon>0$
	\begin{align*}
		J_\rho(u)&\geq\frac{1}{2}\int_{V}|\nabla^su|^2d\mu-\rho C_2(G)\|\nabla^su\|_2-\rho\log C_1(G,h)\\
		&\geq\frac{1}{2}\int_{V}|\nabla^su|^2d\mu-\epsilon\|\nabla^su\|_2^2-\frac{1}{4\epsilon}\le(\rho C_2(G)\ri)^2-\rho\log C_1(G,h).
	\end{align*}
	Taking $\epsilon=1/4$, we arrive at
	\begin{equation}\label{26}
		J_\rho(u)\geq\frac{1}{4}\int_{V}|\nabla^su|^2d\mu-C(G,h,\rho).
	\end{equation}
	
	To sum up, in view of \eqref{25} and \eqref{26}, we conclude that $J_\rho$ has a lower bound in $\mathscr{H}$.\\
	
	\textbf{Step 2.} We shall prove $J_\rho$ has a minimizer $u\in\mathscr{H}$. Take a sequence of function $\{u_k\}\subset\mathscr{H}$ so that
	\begin{equation}\label{27}
		\lim_{k\rightarrow\infty}J_\rho(u_k)=\inf_{u\in\mathscr{H}}J_\rho(u).
	\end{equation}
	For any $\epsilon>0$, by combining \eqref{25}, \eqref{26} and \eqref{27}, we deduce that
	\begin{equation}\label{28}
		\frac{1}{4}\int_{V}|\nabla^su_k|^2d\mu-C(G,h,\rho)\leq J_\rho(u_k)\leq \inf_{u\in\mathscr{H}}J_\rho(u)+\epsilon\leq J_\rho(0)+\epsilon:= C(G,h,\epsilon,\rho),
	\end{equation}
	when $k$ is sufficiently large. Take suitable $\epsilon>0$ such that $C(G,h,\epsilon,\rho)+C(G,h,\rho)>0$. Then it follows from \eqref{28} that for any $\rho\in(-\infty,\,0)\cup(0,\,+\infty)$, there exists a positive constant depending only on $G,h,\rho$ such that
	$$\|\nabla^su_k\|_2^2\leq C(G,h,\rho).$$
	By Lemma \ref{l3}, we have $\{u_k\}$ is bounded in $W^{s,2}(V)$. Hence there exists some $u\in W^{s,2}(V)$ so that, up to a subsequence, $u_k\rightarrow u$ in $W^{s,2}(V)$. Next we only left to verify that $u\in \mathscr{H}$ and $J_\rho(u)=\inf_{u\in\mathscr{H}}J_\rho(u)$. The H\"{o}lder inequality allows 
	\begin{equation}\label{29}
		\le|\int_{V}(u_k-u)d\mu\ri|\leq|V|^\frac{1}{2}\le(\int_{V}|u_k-u|^2d\mu\ri)^\frac{1}{2}\rightarrow0,\quad\text{as}\quad k\rightarrow+\infty,
	\end{equation}
	which yields $\int_{V}ud\mu=0$ and $u\in\mathscr{H}$. On the other hand, by a straightforward computation, one has
	\begin{align}\label{30}
		\nonumber\le|\int_{V}h(e^{u_k}-e^u)d\mu\ri|&\leq\max_{V}h\le|\int_{V}\int_0^1\frac{d}{dt}e^{t(u_k-u)+u}dtd\mu\ri|\\
		\nonumber&=\max_{V}h\le|\int_0^1dt\int_{V}e^{t(u_k-u)+u}(u_k-u)d\mu\ri|\\
		\nonumber&\leq\max_{V}h\int_0^1\le(\int_{V}e^{2t(u_k-u)+2u}d\mu\ri)^\frac{1}{2}\le(\int_{V}|u_k-u|^2d\mu\ri)^\frac{1}{2}dt\\
		\nonumber&\leq\max_{V}h\int_0^1\le(\int_{V}e^{\le(\frac{t(u_k-u)+u}{\|t\nabla^s(u_k-u)+\nabla^su\|_2}\ri)^2+\|t\nabla^s(u_k-u)+\nabla^su\|_2^2}d\mu\ri)^\frac{1}{2}\le(\int_{V}|u_k-u|^2d\mu\ri)^\frac{1}{2}dt\\
		\nonumber&\leq C(G,h)e^{\|\nabla^su_k\|_2^2+\|\nabla^su\|_2^2}\le(\int_{V}|u_k-u|^2d\mu\ri)^\frac{1}{2}\\
		&\rightarrow0,\quad\text{as}\quad k\rightarrow+\infty,
	\end{align}
	where we have used the Cauchy inequality, H\"{o}lder inequality and Moser-Trudinger inequality. Noting that \eqref{29} and \eqref{30}, we obtain that 
	$$J_\rho(u)=\lim_{k\rightarrow\infty}J_\rho(u_k)=\inf_{u\in\mathscr{H}}J_\rho(u).$$
	
	\textbf{Step 3.} Finally, we will show that $u$ satisfies the Euler-Lagrange equation \eqref{15}. Based on variational principle, Lemma \ref{l2} implies that for any $\phi\in\mathscr{H}$
	\begin{align}\label{31}
		\nonumber0&=\dfrac{d}{dt}\bigg{|}_{t=0}J_\rho(u+t\phi)\\
		\nonumber&=\dfrac{d}{dt}\bigg{|}_{t=0}\le(\frac{1}{2}\int_{V}|\nabla^s(u+t\phi)|^2d\mu-\rho\log\int_{V}he^{u+t\phi}d\mu\ri)\\
		&=\int_{V}\le((-\D)^su-\frac{\rho he^u}{\int_{V}he^ud\mu}\ri)\phi d\mu.
	\end{align}
	Let us define the complement space of $\mathscr{H}$ by
	$$\mathscr{H}^\perp=\{u\in W^{s,2}(V):\langle u,\,v\rangle_{W^{s,2}(V)}=0,\ \forall\ v\in\mathscr{H}\}.$$
	Then from \eqref{31} it follows
	\begin{equation}\label{32}
		(-\D)^su-\frac{\rho he^u}{\int_{V}he^ud\mu}\in\mathscr{H}^\perp.
	\end{equation}
	Next we claim that 
	\begin{equation}\label{33}
		\mathscr{H}^\perp=\{\text{const}\}.
	\end{equation}
	It is easy to see that $\{\text{const}\}\subset\mathscr{H}^\perp$. Thus it suffices to prove that $\mathscr{H}^\perp\subset\{\text{const}\}$. For any $f\in\mathscr{H}^\perp$, we denote its integral average over $V$ by $\overline{f}$. Since $(f-\overline{f})\in\mathscr{H}$ and $\overline{f}\in\mathscr{H}^\perp$, it results
	$$0\leq\int_{V}(f-\overline{f})^2d\mu=\int_{V}(f-\overline{f})fd\mu-\int_{V}(f-\overline{f})\overline{f}d\mu=0.$$
	Therefore, we must have $f\equiv\overline{f}$, which leads to $f=\text{const}$, showing that the claim holds. By means of \eqref{32} and \eqref{33}, one can infer that
	\begin{equation}\label{34}
		(-\D)^su=\frac{\rho he^u}{\int_{V}he^ud\mu}+\xi,
	\end{equation}
	where $\xi$ is a constant to be determined. Integrating both sides of \eqref{34} over $V$, utilizing the formula of integration by parts, we find that
	$$\xi=-\frac{\rho}{|V|}.$$
	Hence, we obtain a solution $u\in\mathscr{H}$ satisfying the mean field equation \eqref{15}. This finishes the proof of Theorem \ref{t1}. $\hfill\Box$
	
	\section{The case $0\leq h\not\equiv0$}
	
	In this section, we suppose that $h\in C(V)$ satisfying $0\leq h\not\equiv0$. In order to solve the equation \eqref{15}, let $v=u-\log\int_{V}he^ud\mu$, then \eqref{15} is written in a form
	\begin{equation}\label{38}
		(-\Delta)^s v=\rho he^v-\dfrac{\rho}{|V|},\quad\forall\,x\in V.
	\end{equation}
	Noting that if \eqref{38} has a solution $v$, then \eqref{15} has a solution $u\in\mathscr{M}$, hence we only consider the equation \eqref{38} and look for its solution. 
	
	Now let us introduce the elliptic estimate about the fractional Laplacian $(-\D)^s$, which is a revised version of the corresponding lemma in \cite{L-S-Y}. The following lemma was first proposed by Sun-Wang \cite{S-w} on finite graphs and proved by Li-Sun-Yang \cite{L-S-Y} in the case of $s=1$. For the reader’s convenience, we sketch here. 
	
	\begin{lemma}\label{l8}
		(Elliptic estimate). There exists a positive constant $C$ depending only on the garph $G$ such that for all $u\in C(V)$
		\begin{equation*}
			\max_Vu-\min_Vu\leq C(G)\|(-\Delta)^s u\|_\infty.
		\end{equation*}
	\end{lemma}
	\begin{proof}
		Without loss of generality, we may assume that $V=\{x_1, \cdots,  x_l\}$ and $u(x_1)=\max_{V}u$, $u(x_l)=\min_Vu$. Since for the connectivity of the graph, we may choose $x_1\sim x_2\sim x_3\sim \cdots\sim x_{l-1}\sim x_l$ as the shortest path connecting $x_1$ and $x_l$. By using the Cauchy-Schwartz inequality, we derive that
		\begin{align}\label{41}
			\nonumber0\leq u(x_1)-u(x_l)&\leq\sum_{k=1}^{l-1}|u(x_k)-u(x_{k+1})|\\
			\nonumber&\leq\sqrt{\frac{l-1}{W_s^{0}}}\le(\sum_{k=1}^{l-1}W_s(x_k,x_{k+1})(u(x_k)-u(x_{k+1}))^2\ri)^{\frac{1}{2}}\\
			&\leq\sqrt{\frac{l-1}{W_s^{0}}}\le(\int_{V}|\nabla^s u|^2d\mu\ri)^{\frac{1}{2}}.
		\end{align}
		Here
		$$W_s^{0}=\min_{x,y\in V, x\not=y}W_s(x,y)>0.$$
		From \eqref{06}, we know that
		\begin{equation*}
			\lambda_1^s=\inf_{\overline{v}=0}\frac{\int_{V}|\nabla^s v|^2d\mu}{\int_{V}v^2d\mu}>0
		\end{equation*}
		is the first eigenvalue of $(-\Delta)^s$, where $\lambda_1>0$ is the first eigenvalue of $-\Delta$. It follows from the formula of integration by parts (see Lemma \ref{l2}) and the H\"{o}lder inequality that
		\begin{align*}
			\int_{V}|\nabla^s u|^2d\mu&=\int_{V}(u-\overline{u})(-\D)^s ud\mu\\
			&\leq\le(\int_{V}(u-\overline{u})^2d\mu\ri)^{\frac{1}{2}}\le(\int_{V}((-\D)^su)^2d\mu\ri)^{\frac{1}{2}}\\
			&\leq\le(\frac{1}{\lambda_1^s}\int_{V}|\nabla^s u|^2d\mu\ri)^{\frac{1}{2}}\le(\int_{V}((-\D)^su)^2d\mu\ri)^{\frac{1}{2}}.
		\end{align*}
		This implies that
		\begin{equation}\label{40}
			\int_{V}|\nabla^s u|^2d\mu\leq\frac{1}{\lambda_1^s}\int_{V}((-\D)^su)^2d\mu\leq\frac{|V|}{\lambda_1^s}\|(-\D)^su\|_\infty^2.
		\end{equation}
		Inserting $\eqref{40}$ into $\eqref{41}$, we obtain
		\begin{equation}\label{42}
			\max_Vu-\min_Vu\leq \sqrt{\frac{(l-1)|V|}{W_s^0\lambda_1^s}}\|(-\D)^su\|_\infty,
		\end{equation}
		which ends the proof of Lemma \ref{l8}.
	\end{proof}
	
	The next technical lemma is a priori estimate for solutions of $\eqref{38}$, which is very important for calculation of the topological degree. Our proof depends on the structure of the set $\{x\in V:h(x)=0\}$. Since $h$ is a nonnegative function, there arises new analytical difficulties.
	\begin{lemma}\label{l6}
		Let $G=(V,E)$ be a connected finite graph, $h$ be a function satisfying $0\leq h\not\equiv0$, $s\in(0,\,1)$ and  $\rho\in(-\infty,0)\cup(0,+\infty)$ be some given parameters.  Suppose that there is a real number $\Lambda>0$ such that $h$ and $\rho$ satisfy
		\begin{equation}\label{36}
			|\rho|\leq\Lambda\quad\text{and}\quad\Lambda^{-1}\leq h(x)\leq\Lambda,\quad\forall\ x\in V\setminus\{x\in V:h(x)=0\}.
		\end{equation}
		If $v$ is the solution of $\eqref{38}$, then there exists a constant $C$ depending only on $\Lambda$ and the graph $G$, such that for all $x\in V$, there holds
		$$|v(x)|\leq C(G,\Lambda).$$
	\end{lemma}
	
	\begin{proof}
		We assume that $h$ and $\rho$ satisfy \eqref{36}. If $v$ is a solution of $\eqref{38}$, we consider the cases of $\rho>0$ and $\rho<0$ respectively. \\
		
		\noindent \textbf{Case 1:} $\rho>0$. We first show that $v$ has a uniform upper bound. Since $G$ is a finite graph and $V$ only contains finite vertices, we may assume that there exists $x_1\in V$ so that $v(x_1)=\max_{V}v$. It follows from $\eqref{7}$ and $\eqref{38}$ that
		\begin{equation}\label{39}
			\rho h(x_1)e^{v(x_1)}-\dfrac{\rho}{|V|}=(-\D)^sv(x_1)=\frac{1}{\mu(x_1)}\sum_{y\in V, y\not=x_1}W_s(x_1,y)(v(x_1)-v(y))\geq0.
		\end{equation}
		Noting that $\rho>0$, we deduce from $\eqref{39}$ that $h(x_1)>0$. With this at hand, in view of \eqref{36} and $\eqref{39}$, one obtains 
		\begin{equation*}
			e^{v(x_1)}\geq\frac{1}{|V|h(x_1)}\geq\frac{1}{|V|\Lambda},
		\end{equation*}
		which implies that
		\begin{equation}\label{44}
			\max_{V}v\geq-\log|V|\Lambda.
		\end{equation}
		Integrating both sides of \eqref{38} over $V$, taking into account \eqref{12} and $h(x_1)>0$, we have
		\begin{align*}
			1=\int_{V}he^vd\mu\geq\mu_0h(x_1)e^{v(x_1)},
		\end{align*}
		where $\mu_0=\min_{x\in V}\mu(x)$. This together with \eqref{36} leads to
		\begin{equation}\label{43}
			\max_{V}v\leq\log\frac{\Lambda}{\mu_0}.
		\end{equation}
		
		Secondly, we prove that $v$ also has a uniform lower bound. To see this, by means of \eqref{38}, \eqref{36} and \eqref{43}, one calculates for any $x\in V$
		\begin{align*}
			|(-\D)^sv(x)|&\leq\rho h(x)e^{v(x)}+\frac{\rho}{|V|}\\
			&\leq\Lambda^2e^{\max_{V}v}+\frac{\Lambda}{|V|}\\
			&\leq\frac{\Lambda^3}{\mu_0}+\frac{\Lambda}{|V|}:=a,
		\end{align*}
		which yields 
		\begin{equation*}
			\|(-\D)^sv\|_\infty\leq a.
		\end{equation*}
		From the elliptic estimate \eqref{42} in Lemma \ref{l8}, we know that
		\begin{equation}\label{45}
			\max_Vv-\min_Vv\leq \sqrt{\frac{(l-1)|V|}{W_s^0\lambda_1^s}}a.
		\end{equation}
		Coming back to \eqref{44}, combining with \eqref{45}, we infer that
		\begin{equation}\label{46}
			\min_Vv\geq\max_Vv-\sqrt{\frac{(l-1)|V|}{W_s^0\lambda_1^s}}a\geq-\log|V|\Lambda-\sqrt{\frac{(l-1)|V|}{W_s^0\lambda_1^s}}a.
		\end{equation}
		Hence, in view of \eqref{43} and \eqref{46}, we conclude that
		$$-\log|V|\Lambda-\sqrt{\frac{(l-1)|V|}{W_s^0\lambda_1^s}}a\leq\min_Vv\leq\max_Vv\leq\log\frac{\Lambda}{\mu_0},$$
		as we desired.\\
		
		\noindent \textbf{Case 2:} $\rho<0$. We first present that $v$ is bounded from below. With no loss of generality, we may assume that there exists $x_0\in V$ such that $v(x_0)=\min_{V}v$. Then at $x_0$ it must hold
		\begin{equation*}
			\rho h(x_0)e^{v(x_0)}-\dfrac{\rho}{|V|}=(-\D)^sv(x_0)=\frac{1}{\mu(x_0)}\sum_{y\in V, y\not=x_0}W_s(x_0,y)(v(x_0)-v(y))\leq0.
		\end{equation*}
		Reviewing $\rho<0$, we deduce that $h(x_0)>0$ and in particular,
		\begin{equation}\label{49}
			\min_{V}v\geq-\log|V|\Lambda.
		\end{equation}
		
		Next we prove that $v$ is bounded from above. To this aim, we need to estimate $\max_Vv$. Take $x_1\in V$ so that $v(x_1)=\max_Vv$. Then there are two cases that may occur as follows:\\
		
		\indent\textbf{(i)} $x_1\in\{x\in V:h(x)=0\}$. We claim that $v(x_1)\not=v(x_0)$, precisely, $v$ must not be a constant function. In fact, if $v\equiv\text{const}$, then it follows from \eqref{7} and \eqref{38} that
		$$0=(-\D)^sv(x_1)=-\frac{\rho}{|V|}>0,$$
		which is impossible. Hence at $x_1$ we have
		\begin{align*}
			-\frac{\rho}{|V|}&=(-\D)^sv(x_1)\\
			&=\frac{1}{\mu(x_1)}\sum_{y\in V, y\not=x_1}W_s(x_1,y)(v(x_1)-v(y))\\
			&\geq C(G)\max_{y\in V,y\not=x_1}\le(v(x_1)-v(y)\ri) , 
		\end{align*}
		where we use the fact $\max_{y\in V,y\not=x_1}\le(v(x_1)-v(y)\ri)>0$ since $v\not\equiv\text{const}$, and 
		$$C(G)=\min_{x,y\in V,y\not=x}\frac{W_s(x,y)}{\mu(x)}>0.$$
		Thus it results
		\begin{equation}\label{47}
			0\leq v(x_1)-v(y)\leq-\frac{\rho}{|V|C(G)}\leq\frac{\Lambda}{|V|C(G)},\quad\forall\ y\not=x_1.
		\end{equation}
		Therefore, in order to estimate $v(x_1)$, we take $y=x_0$ in $\eqref{47}$ and have
		\begin{equation}\label{48}
			v(x_1)\leq v(x_0)+\frac{\Lambda}{|V|C(G)}.
		\end{equation}
		Now it is sufficient to establish an upper bound of $v(x_0)$. With the help of $h(x_0)>0$, one can obtain that
		\begin{align*}
			1=\int_{V}he^vd\mu\geq\mu_0h(x_0)e^{v(x_0)},
		\end{align*} 
		which gives
		$$\min_{V}v\leq\log\frac{\Lambda}{\mu_0}.$$
		This together with \eqref{49} and \eqref{48} implies that
		$$-\log|V|\Lambda\leq\min_{V}v<\max_{V}v\leq\log\frac{\Lambda}{\mu_0}+\frac{\Lambda}{|V|C(G)}.$$
		
		\indent\textbf{(ii)} $x_1\in\{x\in V:h(x)>0\}$. It follows from \eqref{38} that
		$$0\leq(-\D)^sv(x_1)=\rho h(x_1)e^{v(x_1)}-\dfrac{\rho}{|V|}.$$
		Under the assumption $\rho<0$, from $h(x_1)>0$, one arrives at 
		\begin{equation*}
			\max_Vv\leq\log\frac{1}{|V|h(x_1)}\leq\log\frac{\Lambda}{|V|}.
		\end{equation*}
		In view of \eqref{49}, we conclude
		$$-\log|V|\Lambda\leq\min_{V}v\leq\max_{V}v\leq\log\frac{\Lambda}{|V|},$$
		which ends the proof of the lemma.
	\end{proof}
	\begin{remark}
		For every positive number $\epsilon$, let $\rho=\pm1$ and $h=\epsilon$ in \eqref{38}, we have
		\begin{equation*}
			(-\D)^s(-\log|V|\epsilon)=0=\pm(\epsilon e^{-\log|V|\epsilon}-\frac{1}{|V|}).
		\end{equation*}
		Thus $v=-\log|V|\epsilon$ is the solution to the equation $(-\D)^sv=\pm(\epsilon e^v-1/|V|)$. However, note that  $$\lim_{\epsilon\rightarrow0^+}-\log|V|\epsilon=+\infty,\quad \lim_{\epsilon\rightarrow+\infty}-\log|V|\epsilon=-\infty.$$ 
		Hence, for any $x\in V\setminus\{x\in V:h(x)=0\}$, the condition $\Lambda^{-1}\leq h(x)\leq\Lambda$ is  necessary.
	\end{remark}
	Let us set $\ell^\infty(V):=\mathscr{X}$. Define a map $\mathscr{F}:\mathscr{X}\rightarrow\mathscr{X}$ by
	\begin{equation}\label{37}
		\mathscr{F}(v)=(-\Delta)^s v-\rho he^v+\frac{\rho}{|V|}.
	\end{equation}
	
	The priori estimate in Lemma \ref{l6} and the homotopy invariance  help us to compute the topological degree of the map $\mathscr{F}$, which is our third lemma in the following.
	\begin{lemma}\label{l7}
		Let $G=(V,E)$ be a finite connected graph, $h$, $s$ and $\rho$ be as in Lemma \ref{l6}, and $\mathscr{F}:\mathscr{X}\rightarrow\mathscr{X}$ be a map given in $\eqref{37}$. Then there exists a large number $R_1>0$ depending only on $\Lambda$ and the graph $G$, such that for all $R\geq R_1$,
		\begin{equation*}
			\deg\le(\mathscr{F},B_R,\theta\ri)= \left\{\begin{aligned}
				&-1, &\rho>0,\\
				&1, &\rho<0,\end{aligned}\right.
		\end{equation*}
		where $B_R=\{u\in\mathscr{X}:\|u\|_\infty<R\}$ is a ball centered at $\theta\in\mathscr{X}$ with radius $R$.		
	\end{lemma}
	\begin{proof}
		Assume that $V=\{x_1,\,\cdots,\,x_n\}$. Then $\mathscr{X}$ is exactly identified with the Euclidean space $\mathbb{R}^n$. We will prove Lemma \ref{l7} with two cases of $\rho>0$ and $\rho<0$ respectively.\\
		
		\noindent \textbf{Case 1:} $\rho>0$. For any given $\epsilon>0$, we define a smooth map $\mathscr{T}_\epsilon:\mathscr{X}\times[0,1]\rightarrow\mathscr{X}$ by
		\begin{equation*}
			\mathscr{T}_\epsilon(v,t)=(-\Delta)^sv-((1-t)\rho+t)((1-t)h+t)e^{v}+\frac{(1-t)\rho+t\epsilon}{|V|},\quad\forall\ (v,t)\in\mathscr{X}\times[0,1].
		\end{equation*}
		Notice that for any $t\in[0,\,1]$
		$$\min\{1,\rho\}\leq(1-t)\rho+t\leq\max\{1,\rho\},\quad\min\{\epsilon,\rho\}\leq(1-t)\rho+t\epsilon\leq\max\{\epsilon,\rho\}.$$
		In addition, for any $(x,t)\in V\setminus\{x\in V:h(x)=0\}\times[0,1]$, there holds 
		$$\min\{\min_{h>0}h,1\}\leq(1-t)h+t\leq\max\{\max_{V}h,1\}.$$
		Hence, there must exist a large number $\Lambda(\epsilon)>0$ depending only on $\epsilon>0$ such that for any $(x,t)\in V\setminus\{x\in V:h(x)=0\}\times[0,1]$
		$$(1-t)\rho+t\leq\Lambda,\quad(1-t)\rho+t\epsilon\leq\Lambda,\quad\Lambda^{-1}\leq (1-t)h(x)+t\leq\Lambda.$$
		Then Lemma \ref{l6} implies that there exists a constant
		$R_\epsilon>1$, depending only on $\Lambda(\epsilon)$ and the graph $G$, such that for all $t\in[0,1]$, all solutions of $\mathscr{T}_\epsilon(v,t)=\theta$ satisfy $\|v_t\|_\infty\leq R_\epsilon$. Denote a ball centered at $\theta\in\mathscr{X}$ with radius $r$ by $B_r\subset\mathscr{X}$, and its boundary by $\p B_r=\{v\in\mathscr{X}:\|v\|_\infty=r\}$. Thus we infer
		$$\p B_R\cap\mathscr{T}_\epsilon^{-1}(\{\theta\})=\varnothing,\quad\forall\ R>R_\epsilon.$$
		Then the topological degree of map $\mathscr{T}_\epsilon$ is well defined for any $R>R_\epsilon$. Now we apply the homotopy invariance of the topological degree to conclude that
		\begin{equation}\label{51}
			\deg\le(\mathscr{T}_\epsilon(\cdotp,0),B_R,\theta\ri)=\deg\le(\mathscr{T}_\epsilon(\cdotp,1),B_R,\theta\ri),\quad\forall\ R>R_\epsilon.
		\end{equation}
		Since $\mathscr{T}_\epsilon(v,0)=\mathscr{F}(v)$, we deduce that
		\begin{equation}\label{52}
			\deg\le(\mathscr{F},B_R,\theta\ri)=\deg\le(\mathscr{T}_\epsilon(\cdotp,0),B_R,\theta\ri),\quad\forall\ R>R_\epsilon.
		\end{equation}
	
		In order to calculate $\deg\le(\mathscr{F},B_R,\theta\ri)$, it is sufficient to calculate $\deg\le(\mathscr{T}_\epsilon(\cdotp,1),B_R,\theta\ri)$. Hence we need to work with the solvability of the equation
		\begin{equation}\label{53}
			\mathscr{T}_\epsilon(v,1)=(-\Delta)^sv-e^v+\frac{\epsilon}{|V|}=\theta.
		\end{equation}
		Next we claim that $v_\epsilon\equiv\log(\epsilon/|V|)$ is the unique solution of \eqref{53} if $\epsilon>0$ is small. Since
		$$\int_Ve^vd\mu=\epsilon,$$
		we have
		\begin{equation}\label{55}
			e^{\max_{V}v_\epsilon}\leq\frac{\epsilon}{\mu_0}.
		\end{equation}
		If $w$ solves 
		$$(-\Delta)^sw=e^w-\frac{\epsilon}{|V|},$$
		then 
		\begin{equation}\label{56}
			(-\Delta)^s(v_\epsilon-w)=e^{v_\epsilon}-e^w,
		\end{equation}
		which leads to
		$$\min_{V}(v_\epsilon-w)<0<\max_{V}(v_\epsilon-w).$$
		As a consequence, we obtain that
		\begin{equation}\label{54}
			|v_\epsilon-w|\leq\max_{V}(v_\epsilon-w)-\min_{V}(v_\epsilon-w).
		\end{equation}
		It follows from \eqref{55}, \eqref{56} and \eqref{54} that for any $x\in V$
		\begin{align}\label{57}
			\nonumber|(-\Delta)^s(v_\epsilon-w)(x)|&=|e^{v_\epsilon}-e^w|\\
			\nonumber&\leq\max\{e^{v_\epsilon},e^w\}|v_\epsilon-w|\\
			&\leq\frac{\epsilon}{\mu_0}\le(\max_{V}(v_\epsilon-w)-\min_{V}(v_\epsilon-w)\ri).
		\end{align}
		With the help of \eqref{42} and \eqref{57}, we have
		\begin{equation}\label{58}
			\max_{V}(v_\epsilon-w)-\min_{V}(v_\epsilon-w)\leq\sqrt{\frac{(l-1)|V|}{W_s^0\lambda_1^s}}\frac{\epsilon}{\mu_0}\le(\max_{V}(v_\epsilon-w)-\min_{V}(v_\epsilon-w)\ri).
		\end{equation} 
		Choose
		$$\epsilon_0=\frac{\mu_0}{2}\sqrt{\frac{W_s^0\lambda_1^s}{(l-1)|V|}}.$$
		If we take $0<\epsilon<\epsilon_0$, then \eqref{58} implies that $v_\epsilon\equiv w$. Therefore, the equation \eqref{53} has a unique solution $v_\epsilon\equiv\log(\epsilon/|V|)$ when $\epsilon\in(0,\,\epsilon_0)$ is small.
		
		Now we rewrite the operator $\mathscr{T}_\epsilon(v,1)$ as follows
		\begin{equation*}
			\mathscr{T}_\epsilon(v,1)=(-\Delta)^s\left(
			\begin{array}{c}
				v(x_1)\\
				\vdots\\
				v(x_n)
			\end{array}\right)-\left(
			\begin{array}{c}
				e^{v(x_1)}\\
				\vdots\\
				e^{v(x_n)}
			\end{array}\right)+\frac{1}{|V|}\left(
			\begin{array}{c}
				\varepsilon\\
				\vdots\\
				\varepsilon
			\end{array}\right)
			.\end{equation*}
		Observing \eqref{06}, we know that the eigenvalues of $(-\D)^s:\mathscr{X}\rightarrow\mathscr{X}$ are represented as
		$$0=\lambda_0^s<\lambda_1^s\leq\lambda_2^s\leq\cdots\leq\lambda^s_{n-1}.$$
		A direct calculation shows
		$$D\mathscr{T}_\epsilon(\log(\epsilon/|V|),1)=\textrm{diag}[0,\, \lambda_1^s,\,\cdots,\,\lambda_{n-1}^s]-\textrm{diag}\le[\frac{\epsilon}{|V|},\,\frac{\epsilon}{|V|},\,\cdots,\,\frac{\epsilon}{|V|}\ri],$$
		where we denote the $n\times n$ diagonal matrix by $\textrm{diag}[\cdot,\,\cdot,\,\ldots,\,\cdot]$. We may choose a small
		$$\epsilon_1=\min\{\epsilon_0,\lambda_1^s|V|\},$$
		such that for any $\epsilon\in(0,\epsilon_1)$, $\mathscr{T}_\epsilon(v,1)=\theta$ has a unique solution $v_\epsilon$ satisfying $e^{v_\epsilon}=\epsilon/|V|<\lambda_1^s$. Hence, for some given $\epsilon\in(0,\epsilon_1)$, we deduce that
		$$\deg\le(\mathscr{T}_\epsilon(\cdotp,1),B_R,\theta\ri)=\textrm{sgn}\det{D\mathscr{T}_\epsilon(\log(\epsilon/|V|),1)}=\textrm{sgn}\le\{-\frac{\epsilon}{|V|}\prod_{j=1}^{n-1}(-\frac{\epsilon}{|V|}+\lambda^s_j)\ri\}=-1.$$
		Finally, it follows from $\eqref{51}$ and $\eqref{52}$ that
		\begin{align*}
			\deg\le(\mathscr{F},B_R,\theta\ri)=-1,\quad\forall\ R>R_\epsilon.
		\end{align*}
		
		\noindent \textbf{Case 2:} $\rho<0$. Define a smooth map $\mathscr{T}:\mathscr{X}\times[0,1]\rightarrow\mathscr{X}$ by
		\begin{equation}\label{50}
			\mathscr{T}(v,t)=(-\Delta)^sv-((1-t)\rho-t)((1-t)h+t)e^{v}+\frac{(1-t)\rho-t}{|V|},\quad\forall\ (v,t)\in\mathscr{X}\times[0,1].
		\end{equation}
		From the above argument, there exists a large number $\Lambda>0$ such that 
		$$-\Lambda\leq(1-t)\rho-t,\quad\Lambda^{-1}\leq (1-t)h(x)+t\leq\Lambda,\quad\forall\ (x,t)\in V\setminus\{x\in V:h(x)=0\}\times[0,1].$$
		Applying Lemma \ref{l6} again, we can also find a number $R_0>1$ depending only on $\Lambda$ and the graph $G$ such that all solutions $v_t$ of $\mathscr{T}(v,t)=\theta$ is uniformly bounded for any $t\in[0,\,1]$, where $\mathscr{T}(v,t)$ is given in \eqref{50}. Thus the topological degree of map $\mathscr{T}$ is well defined for any $R\geq R_0$. By using the homotopy invariance, we have
		\begin{equation}\label{59}
			\deg\le(\mathscr{F},B_R,\theta\ri)=\deg\le(\mathscr{T}(\cdotp,0),B_R,\theta\ri)=\deg\le(\mathscr{T}(\cdotp,1),B_R,\theta\ri),\quad\forall\ R>R_0,
		\end{equation}
		where $\mathscr{F}$ is defined by \eqref{37}.
		
		Next we shall compute $\deg\le(\mathscr{T}(\cdotp,1),B_R,\theta\ri)$. We only need to treat the following equation
		\begin{equation}\label{60}
			(-\D)^s v=-e^{v}+\frac{1}{|V|}.
		\end{equation}
		We claim that $v\equiv-\log|V|$ is the unique solution of $\eqref{60}$. Indeed, $v$ cannot be anything but constant function. Suppose not, we can look for some $x_0$, $x_1\in V$ such that $v(x_0)=\min_Vv$, $v(x_1)=\max_Vv$ satisfying $v(x_0)<v(x_1)$. From \eqref{7}, one can easily see that
		\begin{equation*}
			(-\D)^sv(x_0)=\frac{1}{\mu(x_0)}\sum_{y\in V, y\not=x_0}W_s(x_0,y)(v(x_0)-v(y))<0,
		\end{equation*}
		and
		\begin{equation*}
			(-\D)^sv(x_1)=\frac{1}{\mu(x_1)}\sum_{y\in V, y\not=x_1}W_s(x_1,y)(v(x_1)-v(y))>0.
		\end{equation*}
		Thus we have
		$$-e^{\min_Vv}+\frac{1}{|V|}=(-\D)^sv(x_0)<0,$$
		and
		$$-e^{\max_Vv}+\frac{1}{|V|}=(-\D)^sv(x_1)>0,$$
		which yields
		$$-e^{\min_Vv}<-\frac{1}{|V|}<-e^{\max_Vv}.$$
		However, this is impossible, hence confirming our claim.
		
		Now write the operator $\mathscr{T}(v,1):\mathscr{X}\rightarrow\mathscr{X}$ as 
		\begin{equation*}
			\mathscr{T}(v,1)=(-\Delta)^s\left(
			\begin{array}{c}
				v(x_1)\\
				\vdots\\
				v(x_n)
			\end{array}\right)+\left(
			\begin{array}{c}
				e^{v(x_1)}\\
				\vdots\\
				e^{v(x_n)}
			\end{array}\right)-\frac{1}{|V|}\left(
			\begin{array}{c}
				1\\
				\vdots\\
				1
			\end{array}\right)
			.\end{equation*}
		Therefore, a direct computation derives that
		$$D\mathscr{T}(-\log|V|,1)=\textrm{diag}[0,\, \lambda_1^s,\,\cdots,\,\lambda_{n-1}^s]+\textrm{diag}\le[\frac{1}{|V|},\,\frac{1}{|V|},\,\cdots,\,\frac{1}{|V|}\ri].$$
		It is obvious that
		$$\deg\le(\mathscr{T}(\cdotp,1),B_R,\theta\ri)=\textrm{sgn}\det{D\mathscr{T}(-\log|V|,1)}=\textrm{sgn}\le\{\frac{1}{|V|}\prod_{j=1}^{n-1}(\frac{1}{|V|}+\lambda^s_j)\ri\}=1,\quad\forall\ R>R_0.$$
		Noting that \eqref{59}, we have
		\begin{align*}
			\deg\le(\mathscr{F},B_R,\theta\ri)=1,\quad\forall\ R>R_0.
		\end{align*}
		
		To sum up, for fixed $\epsilon\in(0,\epsilon_1)$, taking $R_1=\max\{R_\epsilon,R_0\}$, we complete the proof of the lemma.
		
	\end{proof}

	As an application of the topological degree presented in Lemma \ref{l7}, for the case $h$ is nonnegative, we give the answer to problem concerning the existence results of $\eqref{15}$.\\
	
	\noindent$\textbf{\emph{Proof of Theorem \ref{t2}}.}$ For any $\rho\in(-\infty,\,0)\cup(0,\,+\infty)$ and $s\in(0,\,1)$, we can use Lemma \ref{l7} to find some large $R_1>1$ so that
	$$\deg\le(\mathscr{F},B_{R_1},\theta\ri)\not=0.$$
	Thus the Kronecker’s existence theorem implies that $\eqref{38}$ has at least one solution $v$. Hence, \eqref{15} has at least one solution $u\in\mathscr{M}$ satisfying $\int_{V}he^ud\mu=1$. This ends the proof of Theorem \ref{t2}. $\hfill\Box$

	\section{The case $\max_Vh>0$}
	
	In this section, we suppose that $h$ is a sign-changing prescribed function, and prove Theorems \ref{t3} and \ref{t4} in two subsections.

	\subsection{The first flow}
	
	In the spirit of our recent work \cite{L-Z}, we want to revisit the existence of the solutions for \eqref{15}, and no longer require that $h$ is nonnegative. Recall that the first mean field type heat flow is given by
	\begin{equation}\label{61}\left\{\begin{array}{lll}
			\dfrac{\p}{\p t}e^u=-(-\Delta)^s u+\rho\left(\dfrac{he^u}{\int_V he^ud\mu}-\dfrac{1}{|V|}\right),\quad\left(x,\, t\right)\,\in \, V\times\left(0,\, +\infty\right),\\[3ex]
			u\left(x,\, 0\right)=u_0\left(x\right),\quad\quad\quad\quad\quad\quad\quad\ \  x\,\in \,V,\end{array}\right.
	\end{equation}
	where $\rho>0$. It is easy to check that the first flow is a negative gradient flow of the functional
	\begin{equation}\label{62}
		J_{\rho,h}\left(u\right)=\frac{1}{2}\int_V |\nabla^s u|^2d\mu-\rho\log\bigg{|}\int_V he^ud\mu\bigg{|}+\frac{\rho}{|V|}\int_V ud\mu,\quad\forall\ u\in W^{s,2}\left(V\right).
	\end{equation}
	In particular, (\ref{62}) is the variational functional associated with the mean field equation \eqref{15}.

	\subsubsection{Short time existence}   
	
	Firstly, we focus on the short time existence of the first flow \eqref{61}. We assume $V=\{x_1,\, x_2,\, \cdots,\, x_n\}$ for some integer $n\geq1$. Then any function $u\in C(V)$ can be characterized as $\textbf{u}=\left(u\left(x_1\right),\ \cdots,\ u\left(x_n\right)\right)\in\mathbb{R}^n$.
	Let us denote
	\begin{equation}\label{66}
		F\left(u\right)(x)=-(-\Delta)^s u(x)+\rho\left(\dfrac{he^{u(x)}}{\int_V he^ud\mu}-\dfrac{1}{|V|}\right),
	\end{equation}
	and define a map $\textbf{T}:\mathbb{R}^n\rightarrow\mathbb{R}^n$ by $\textbf{T}\left(\textbf{u}\right)=\left(\textbf{T}_1\left(\textbf{u}\right),\, \cdots,\, \textbf{T}_n\left(\textbf{u}\right)\right)$, where $\textbf{T}_j\left(\textbf{u}\right)=F(u)(x_j)$ for every $1\leq j\leq n$. Clearly, one has
	$$\textbf{T}\left(\textbf{u}\right)=\le(-(-\Delta)^s u(x_1)+\rho\left(\dfrac{he^{u(x_1)}}{\int_V he^ud\mu}-\dfrac{1}{|V|}\right),\,\cdots,\,-(-\Delta)^s u(x_n)+\rho\left(\dfrac{he^{u(x_n)}}{\int_V he^ud\mu}-\dfrac{1}{|V|}\right)\ri)\in\mathbb{R}^n.$$
	Then the flow \eqref{61} is equivalent to the ordinary differential system
	\begin{equation}\label{67}\left\{\begin{array}{lll}
			\dfrac{d}{dt}e^{u\left(x_1\right)}=\textbf{T}_1\left(\textbf{u}\right),\\[1ex]
			\vdots\\[1ex]
			\dfrac{d}{dt}e^{u\left(x_n\right)}=\textbf{T}_n\left(\textbf{u}\right),\\[1ex]
			\textbf{u}\left(0\right)=\textbf{u}_0.\end{array}\right.
	\end{equation}
	Here $\textbf{u}_0=\left(u_0\left(x_1\right),\, \cdots,\, u_0\left(x_n\right)\right)$ is the initial value. Now we have
	
	\begin{lemma}\label{l9}
		Let $G=\left(V,\, E\right)$ be a connected finite graph. For any initial function $u_0\in C(V)$, there exists some constant $T_0>0$ such that $u:V\times[0,\, T_0]\rightarrow\mathbb{R}$ is a unique solution of the first flow (\ref{61}).
	\end{lemma}
	\begin{proof}
		We claim that $\textbf{T}:\mathbb{R}^n\rightarrow\mathbb{R}^n$ is analytic. For every $1\leq j\leq n$, $\textbf{T}_j\left(\textbf{u}\right)$ can be written as
		\begin{equation}\label{143}
			\textbf{T}_j\left(\textbf{u}\right)=F(u)(x_j)=\frac{1}{\mu(x_j)}\sum_{y\in V, y\not=x_j}W_s(x_j,y)(u(y)-u(x_j))+\rho\left(\dfrac{h(x_j )e^{u (x_j )}}{\sum_{i=1}^n\mu\left(x_i\right)h\left(x_i\right)e^{u\left(x_i\right)}}-\dfrac{1}{|V|}\right).
		\end{equation}
		Observe that $\left(u\left(x_1\right),\, \cdots,\, u\left(x_n\right)\right)$ on the right hand of \eqref{143} can be replaced by $\textbf{u}$. This shows that $\textbf{T}_j$ is analytic. Then the short time existence theorem of ordinary differential equation implies that there exist some $T_0>0$ and a $C^1$ map $\textbf{u}:[0,\, T_0]\rightarrow\mathbb{R}^n$ such that $\textbf{u}\left(t\right)$ is a unique solution of (\ref{67}) in a time interval $[0,\, T_0]$. Let us set $\textbf{u}(t)=(u(x_1,\,  t),\, \cdots,\, u(x_n,\, t))$, then $u:V\times[0,\, T_0]\rightarrow\mathbb{R}$ is a unique solution of the first heat flow (\ref{61}), which ends the proof of Lemma \ref{l9}.
	\end{proof}
	
	\subsubsection{Long time existence}   
	
	Once Lemma \ref{l9} holds, we suppose that
	\begin{equation}\label{68}
		T=\sup\Big{\{}T_0>0\ |\ u:V\times[0,\, T_0]\rightarrow\mathbb{R}\ \mathrm{is\ the\ unique\ solution\ of\ the\ flow\ \eqref{61}}\Big{\}}<+\infty.
	\end{equation}
	Now we would like to give two important properties of the flow \eqref{61} for any $t\in[0,\, T)$, i.e.
	\begin{proposition}\label{p4}
		Let $T$ be defined in (\ref{68}). Then we have the following properties.\\
		(i) (Invariant quality) For any $t\in[0,\, T)$, there holds
		\begin{equation}\label{69}
			\int_V e^{u\left(\cdot,\, t\right)}d\mu=\int_V e^{u_0}d\mu.
		\end{equation}
		(ii) (Monotone) For any $0\leq t_1<t_2<T$, it must hold
		\begin{equation}\label{70}
			J_{\rho,h}\left(u\left(\cdot,\,t_2\right)\right)\leq J_{\rho,h}\left(u\left(\cdot,\,t_1\right)\right).
		\end{equation}
	\end{proposition}
	\begin{proof}
		$(i)$ For any $\left(x,\, t\right)\in V\times[0,\, T)$, it follows from the definition of $T$ in \eqref{68} that $u\left(x,\, t\right)$ is a solution of \eqref{61}. Differentiating with respect to $t$, for any $t\in[0,\, T)$, we infer that
		\begin{equation*}
			\begin{split}
				\dfrac{d}{dt}\int_V e^{u\left(\cdot,\, t\right)}d\mu&=\int_V\dfrac{\p}{\p t}e^{u\left(\cdot,\, t\right)}d\mu\\
				&=-\int_V(-\Delta^s) u\left(\cdot,\, t\right)d\mu+\int_V\rho\left(\dfrac{he^{u\left(\cdot,\, t\right)}}{\int_V he^{u\left(\cdot,\, t\right)}d\mu}-\dfrac{1}{|V|}\right)d\mu\\
				&=0,
			\end{split}
		\end{equation*}
		where we have used the formula of integration by parts \eqref{12}. Hence we obtain 
		\begin{equation*}
			\int_V e^{u\left(\cdot,\, t\right)}d\mu=\int_V e^{u\left(\cdot,\, 0\right)}d\mu=\int_V e^{u_0}d\mu.
		\end{equation*}
		
		\noindent$(ii)$ In view of (\ref{62}) and (\ref{66}), we deduce that
		\begin{equation}\label{71}
			\dfrac{d}{dt}J_{\rho,h}\left(u\left(\cdot,\, t\right)\right)=-\int_VF\left(u\right)\dfrac{\p u}{\p t} d\mu=-\int_Ve^{u\left(\cdot,\, t\right)}\bigg{|}\dfrac{\p u}{\p t}\bigg{|}^2 d\mu\leq0.
		\end{equation}
		Thus it results that \eqref{70} holds.
	\end{proof}
	
	For all $t$ whenever the flow \eqref{61} exists, we state and prove the crucial compactness result of the solution $u\left(\cdot,\,t\right)$, which allows us to establish the long time existence of \eqref{61}. 
	
	\begin{lemma}\label{l10}
		For any $t\in[0,\,T)$, suppose that $u\left(x,\, t\right)$ is the solution of the first flow. Then there exists a constant $C$ depending only on $G$, $T$, $\rho$, $h$ and $u_0$ such that
		\begin{equation}\label{72}
			\|u\left(\cdot,\,t\right)\|_{W^{s,2}\left(V\right)}\leq C\left(G,\, T,\, \rho,\, h,\, u_0\right),\quad\forall\, t\in[0,\,T),
		\end{equation}
		where the norm $\|\cdot\|_{W^{s,2}\left(V\right)}$ is defined by (\ref{73}).
	\end{lemma}
	\begin{proof}
		Let $T$ be given as (\ref{68}), and $T<+\infty$ be a finite real number. We seperate the proof into four steps.\\
		
		\textbf{Step 1.} As a first step, we claim that there exist some constants such that for any $t\in[0,\, T)$,
		\begin{equation}\label{74}
			0<C\left(G,\, \rho,\, h,\, u_0\right)\leq\left|\int_V he^{u\left(\cdot,\, t\right)}d\mu\right|\leq C\left(G,\, h,\, u_0\right).
		\end{equation}
		Using Lemma $\ref{l4}$, for any $\beta>0$ and $t\in[0,\, T)$,  there exists a constant $C(G,\beta)$ such that
		\begin{equation}\label{75}
			\int_V e^{\beta\frac{\left(u\left(\cdot,\, t\right)-\overline{u}\left(t\right)\right)^2}{\|\nabla^s u\left(\cdot,\, t\right)\|^2_{2}}}d\mu\leq C\left(G,\, \beta\right),
		\end{equation}
		where $\bar{u}\left(t\right)= \int_{V}u\left(\cdot,\, t\right)d\mu/{|V|}$. It follows from (\ref{75}) and the Young inequality that for any $\epsilon>0$ and $t\in[0,\, T)$
		\begin{align*}
			\log\int_Ve^{u\left(\cdot,\, t\right)-\overline{u}\left(t\right)}d\mu&
			\leq\log\int_V e^{\frac{\left(u\left(\cdot,\, t\right)-\overline{u}\left(t\right)\right)^2}{4\epsilon\|\nabla^s u\left(\cdot,\, t\right)\|^2_2}+\epsilon\|\nabla^s u\left(\cdot,\ t\right)\|_{2}^2}d\mu\\
			&\leq\epsilon\int_{V}|\nabla^s u\left(\cdot,\, t\right)|^2d\mu+C\left(G,\, \epsilon\right).
		\end{align*}
		This together with $\rho>0$ leads to
		\begin{equation}\label{76}
			\frac{\rho}{|V|}\int_{V}u\left(\cdot,\, t\right)d\mu\geq\rho\log\int_Ve^{u\left(\cdot,\, t\right)}d\mu-\epsilon\rho\int_{V}|\nabla^s u\left(\cdot,\, t\right)|^2d\mu+C\left(G,\, \epsilon,\, \rho\right).
		\end{equation}
		Inserting (\ref{76}) into (\ref{62}), noting that (\ref{70}), we derive that for any $t\in[0,\, T)$,
		\begin{align}\label{77}
			\nonumber	J_{\rho,h}\left(u_0\right) \geq& J_{\rho,h}\left(u\left(\cdot,\, t\right)\right)\\
			\geq&\left(\frac{1}{2}-\epsilon\rho\right)\int_{V}|\nabla^s u\left(\cdot,\, t\right)|^2d\mu-\rho\log\left|\int_V he^{u\left(\cdot,\, t\right)}d\mu\right|+\rho\log\int_Ve^{u\left(\cdot,\, t\right)}d\mu+C\left(G,\, \epsilon,\, \rho\right).
		\end{align}
		Let us take $\epsilon=1/\left(2\rho\right)$ in (\ref{77}). Then from (\ref{69}), one has
		\begin{equation}\label{93}
			J_{\rho,h}\left(u_0\right)\geq-\rho\log\left|\int_V he^{u\left(\cdot,\, t\right)}d\mu\right|+\rho\log\int_Ve^{u_0}d\mu+C\left(G,\, \rho\right),
		\end{equation}
		which implies that
		$$\log\left|\int_V he^{u\left(\cdot,\, t\right)}d\mu\right|\geq\log\int_Ve^{u_0}d\mu-\frac{1}{\rho}J_{\rho,h}\left(u_0\right)+C\left(G,\, \rho\right).$$
		Hence, we infer that for any $t\in[0,\, T)$,
		\begin{equation}\label{144}
			\left|\int_V he^{u\left(\cdot,\, t\right)}d\mu\right|\geq e^{-\frac{1}{\rho}J_{\rho,h}\left(u_0\right)+C\left(G,\, \rho\right)}\int_Ve^{u_0}d\mu:=C\left(G,\, \rho,\, h,\, u_0\right)>0.
		\end{equation}
		On the other hand, it follows from (\ref{69}) that for any $t\in[0,\, T)$,
		\begin{equation}\label{145}
			\left|\int_V he^{u\left(\cdot,\, t\right)}d\mu\right|\leq\max_{V}|h|\int_Ve^{u_0}d\mu:= C\left(G,\, h,\, u_0\right).
		\end{equation}
		Combining \eqref{144} and \eqref{145}, \eqref{74} follows immediately.
		
		\textbf{Step 2.} Now we study the linearization of the first flow of the form
		\begin{equation*}\left\{\begin{array}{lll}
				\dfrac{\p u}{\p t}=e^{-u}F\left(u\right),\quad\left(x,\, t\right)\,{\in}\, V\times[0,\, T),\\[2ex]
				u\left(x,\, 0\right)=u_0\left(x\right),\quad x\,{\in}\,V,\end{array}\right.
		\end{equation*}
		where $F$ is defined by \eqref{66}. By applying Lemma \ref{l2} and Proposition \ref{p4}, we get
		\begin{align}\label{78}
			\nonumber\frac{d}{dt}\int_V e^{2u\left(\cdot,\, t\right)}d\mu&=2\int_V e^{u\left(\cdot,\, t\right)}F\left(u\right)d\mu\\
			&=-2\int_V \nabla^se^{u\left(\cdot,\, t\right)}\nabla^su\left(\cdot,\, t\right)d\mu+2\rho\left(\dfrac{\int_Vhe^{2u\left(\cdot,\, t\right)}d\mu}{\int_V he^{u\left(\cdot,\, t\right)}d\mu}-\dfrac{1}{|V|}\int_Ve^{u_0}d\mu\right).
		\end{align}
		Noting that \eqref{13}, one knows that
		\begin{equation*}
			\int_V \nabla^se^{u\left(\cdot,\, t\right)}\nabla^su\left(\cdot,\, t\right)d\mu=\frac{1}{2}\sum_{x\in V}\sum_{y\in V, y\not=x}W_s(x,y)\left(e^{u\left(x,\, t\right)}-e^{u\left(y,\, t\right)}\right)\left(u\left(x,\, t\right)-u\left(y,\, t\right)\right).
		\end{equation*}	
		It yields that
		\begin{equation}\label{79}
			\int_V \nabla^se^{u\left(\cdot,\, t\right)}\nabla^su\left(\cdot,\, t\right)d\mu\geq0,\quad\forall\, t\in[0,\, T).
		\end{equation}
		Combining (\ref{78}) and (\ref{79}), in view of (\ref{74}), we arrive at
		\begin{equation}\label{80}
			\frac{d}{dt}\int_V e^{2u\left(\cdot,\, t\right)}d\mu\leq C\left(G,\, \rho,\, h,\, u_0\right)\int_V e^{2u\left(\cdot,\, t\right)}d\mu.
		\end{equation}
		Integrating (\ref{80}) from $0$ to $t$, by a simple calculation, we deduce that
		\begin{equation}\label{81}
			\int_V e^{2u\left(\cdot,\, t\right)}d\mu\leq e^{Ct}\int_V e^{2u_0}d\mu\leq C\left(G,\, T,\, \rho,\, h,\, u_0\right),\quad\forall\ t\in[0,\, T).
		\end{equation}
		
		\textbf{Step 3.} Set
		\begin{equation}\label{146}
			\int_Ve^{u_0}d\mu=c_0>0.
		\end{equation}
		For any $\epsilon>0$ and $t\in[0,\, T)$, let us define a set
		\begin{equation}\label{82}
			V_{\epsilon,t}=\left\{x\in V\ |\ e^{u\left(\cdot,\, t\right)}\geq\epsilon\right\}.
		\end{equation}
		With the help of (\ref{79}), \eqref{81} and (\ref{82}), by using the H\"{o}lder inequality, one derives that
		\begin{align}\label{83}
			\nonumber		c_0&=\int_V e^{u\left(\cdot,\, t\right)}d\mu\\
			\nonumber	 &=\int_{V_{\epsilon,t}} e^{u\left(\cdot,\, t\right)}d\mu+\int_{V\setminus V_{\epsilon,t}} e^{u\left(\cdot,\, t\right)}d\mu\\
			\nonumber			&\leq\epsilon|V\setminus V_{\epsilon,t}|+\left(\int_{V_{\epsilon,t}} e^{2u\left(\cdot,\, t\right)}d\mu\right)^{\frac{1}{2}}\left(\int_{V_{\epsilon,t}} 1d\mu\right)^{\frac{1}{2}}\\
			&\leq\epsilon|V|+C\left(G,\, T,\, \rho,\, h,\, u_0\right)|V_{\epsilon,t}|^{\frac{1}{2}}.
		\end{align}
		Taking $\epsilon=c_0/({2|V|})$ in (\ref{83}), we have
		\begin{equation}\label{84}
			|V_{c_0/({2|V|}),t}|\geq C\left(G,\, T,\, \rho,\, h,\, u_0\right)>0.
		\end{equation}
		With this at hand, we obtain $V_{c_0/({2|V|}),t}\not=\varnothing$. From now on, denote $A_t=V_{c_0/({2|V|}),t}$. By means of (\ref{82}) and  (\ref{84}), we infer that
		\begin{equation}\label{91}
			\int_{A_t} u\left(\cdot,\, t\right)d\mu\geq C\left(G,\, T,\, \rho,\, h,\, u_0\right)\log\dfrac{c_0}{2|V|}.
		\end{equation}
		Now we estimate the upper bound of $\int_{A_t} u\left(\cdot,\, t\right)d\mu$. Using the Jensen’s inequality and the H\"{o}lder inequality, from \eqref{81} and \eqref{84} we have
		\begin{align*}
			e^{\frac{\int_{A_t} u\left(\cdot,\, t\right)d\mu}{|A_t|}}=e^{\sum_{x\in A_t}\frac{\mu(x)}{|A_t|}u\left(\cdot,\, t\right)}\leq\sum_{x\in A_t}\frac{\mu(x)}{|A_t|}e^{u\left(\cdot,\, t\right)}&=\frac{1}{|A_t|}\int_{A_t}e^{u\left(\cdot,\, t\right)}d\mu\\
			&\leq\frac{1}{|A_t|}\le(\int_{A_t}e^{2u\left(\cdot,\, t\right)}d\mu\ri)^{\frac{1}{2}}\le(\int_{A_t}1d\mu\ri)^{\frac{1}{2}}\\
			&\leq\frac{1}{|A_t|^\frac{1}{2}}\le(\int_{V}e^{2u\left(\cdot,\, t\right)}d\mu\ri)^{\frac{1}{2}}\leq C\left(G,\, T,\, \rho,\, h,\, u_0\right).
		\end{align*}
		Consequently, it must hold
		\begin{equation}\label{92}
			\int_{A_t} u\left(\cdot,\, t\right)d\mu\leq|V|\log C\left(G,\, T,\, \rho,\, h,\, u_0\right).
		\end{equation}
		Combining \eqref{91} and \eqref{92}, we conclude that
		\begin{equation}\label{99}
			\bigg{|}\int_{A_t} u\left(\cdot,\, t\right)d\mu\bigg{|}\leq C\left(G,\, T,\, \rho,\, h,\,u_0\right),
		\end{equation}
		which together with the Young inequality and the H\"{o}lder inequality gives that for any $\epsilon>0$
		\begin{align}\label{85}
			\nonumber\left(\int_V u\left(\cdot,\, t\right)d\mu \right)^2
			&=\left(\int_{A_t} u\left(\cdot,\, t\right)d\mu\right)^2+\left(\int_{V\setminus A_t} u\left(\cdot,\, t\right)d\mu\right)^2+2\int_{A_t} u\left(\cdot,\, t\right)d\mu\int_{V\setminus A_t} u\left(\cdot,\, t\right)d\mu\\
			\nonumber	&\leq\left(\int_{A_t} u\left(\cdot,\, t\right)d\mu\right)^2+\left(\int_{V\setminus A_t} u\left(\cdot,\, t\right)d\mu\right)^2\\
			\nonumber&\quad+\frac{1}{\epsilon}\left(\int_{A_t} u\left(\cdot,\, t\right)d\mu\right)^2
			+\epsilon\left(\int_{V\setminus A_t} u\left(\cdot,\,    t\right)d\mu\right)^2\\
			&\leq\left(1+\frac{1}{\epsilon}\right)C\left(G,\, T,\, \rho,\, h,\, u_0\right)+\left(1+\epsilon\right)|V\setminus A_t|\int_V u^2\left(\cdot,\, t\right)d\mu.
		\end{align}	
		It follows from Lemma \ref{l3} that there exists a constant $C(G)$ such that
		$$\int_V \left(u\left(\cdot,\, t\right)-\overline{u}\left(t\right)\right)^2d\mu\leq C(G)\int_V |\nabla^s u\left(\cdot,\, t\right)|^2d\mu.$$
		This together with (\ref{85}) implies that
		\begin{align}\label{86}
			\nonumber\int_V u^2\left(\cdot,\, t\right)d\mu\leq& C(G)\int_V |\nabla^s u\left(\cdot,\, t\right)|^2d\mu+\int_V \overline{u}^2\left(\cdot,\, t\right)d\mu\\
			\leq& C(G)\int_V |\nabla^s u\left(\cdot,\, t\right)|^2d\mu+\dfrac{\left(1+\epsilon\right)|V\setminus A_t|}{|V|}\int_V u^2\left(\cdot,\, t\right)d\mu+C(G,\, T,\, \rho,\, h,\, u_0,\,\epsilon).
		\end{align}
		By taking suitable $\epsilon>0$ in (\ref{86}), we conclude that there exists some constant $C$ depending only on $G$, $T$, $\rho$, $h$ and $u_0$ such that
		\begin{equation}\label{87}
			\int_V u^2\left(\cdot,\, t\right)d\mu\leq C\left(G\right)\int_{V}|\nabla^s u\left(\cdot,\, t\right)|^2d\mu+C\left(G,\, T,\, \rho,\, h,\, u_0\right),\quad\forall\, t\in[0,\, T).
		\end{equation}
		
		\textbf{Step 4.} According to \eqref{79} and \eqref{80}, we obtain from \eqref{62} that for any $t\in[0,\, T)$
		\begin{align}\label{88}
			\nonumber\frac{1}{2}\int_V |\nabla^s u\left(\cdot,\, t\right)|^2d\mu+\frac{\rho}{|V|}\int_V u\left(\cdot,\, t\right)d\mu&\leq J_{\rho,h}\left(u_0\right)+\rho\log\left|\int_V he^{u\left(\cdot,\, t\right)}d\mu\right|\\
			\nonumber&\leq J_{\rho,h}\left(u_0\right)+\rho\left(\log\max_{V}|h|+\log\int_V e^{u_0}d\mu\right)\\
			&\leq C\left(G,\, \rho,\, h,\, u_0\right).
		\end{align}
		Applying the Young inequality to \eqref{88}, one has for any $\epsilon>0$
		\begin{align}\label{89}
			\nonumber\int_V |\nabla^s u\left(\cdot,\, t\right)|^2d\mu&\leq C\left(G,\, \rho,\, h,\, u_0\right)-\frac{2\rho}{|V|}\int_V u\left(\cdot,\, t\right)d\mu\\
			&\leq C\left(G,\, \rho,\, h,\, u_0\right)+\frac{2\rho\epsilon}{|V|}\int_V u^2\left(\cdot,\, t\right)d\mu+\frac{\rho}{2\epsilon}.
		\end{align}
		Inserting (\ref{87}) into (\ref{89}), we get
		$$\int_V |\nabla^s u\left(\cdot,\, t\right)|^2d\mu-\frac{\rho}{2\epsilon}-C\leq\frac{2\rho\epsilon}{|V|}\int_V u^2\left(\cdot,\, t\right)d\mu\leq\frac{2\rho C\epsilon}{|V|}\left(\int_{V}|\nabla^s u\left(\cdot,\, t\right)|^2d\mu+1\right).$$
		Taking $\epsilon=|V|/\left(4\rho C\right)$, we deduce that
		\begin{equation}\label{90}
			\int_V |\nabla^s u\left(\cdot,\, t\right)|^2d\mu\leq C\left(G,\, T,\, \rho,\, h,\, u_0\right).
		\end{equation}
		Combining \eqref{87} and (\ref{90}), we conclude that
		$$\|u\left(\cdot,\,t\right)\|_{W^{s,2}\left(V\right)}\leq C\left(G,\, T,\, \rho,\, h,\, u_0\right),\quad\forall\ t\in[0,\,T),$$
		which completely ends the proof of lemma.	
	\end{proof}
	With the help of Lemma \ref{l10}, now we are able to give the long time existence of the first flow.\\
	
	\noindent$\textbf{\emph{Proof of (i) in Theorem \ref{t3}}.}$ Suppose that $T<+\infty$. Then there exists a unique solution $u(x,\,t)$ of the first flow in a time interval $[0,\,T)$. By employing Lemma \ref{l10} and the short time existence theorem of ordinary differential equation, $u\left(x,\, t\right)$ can be uniquely extended beyond $T$, in contradiction with the deifinition of $T$ in (\ref{68}). Therefore, $T=+\infty$. Hence $(i)$ of Theorem \ref{t3} is proved. $\hfill\Box$
	
	\subsubsection{Global convergence}
	
	In this part, our aim is to obtain the convergence results of the solution $u:V\times[0,\, +\infty)\rightarrow\mathbb{R}$ for the first flow \eqref{61}. Our approach can be divided into two aspects. One is to find an increasing sequence $t_n\rightarrow+\infty$ and some function $u_\infty\in C(V)$ such that $u\left(x,\, t_n\right)\rightarrow u_\infty\left(x\right)$ uniformly in $x\in V$ as $n\rightarrow+\infty$. The other is to prove that along the flow (\ref{61}), $u\left(x,\, t\right)\rightarrow u_\infty\left(x\right)$ as $t\rightarrow+\infty$ uniformly on $V$.
	
	Observing that $J_\rho\left(u\left(\cdot,\, t\right)\right)$ decreases along the flow, integrating (\ref{71}) from $0$ to $t$, making use of \eqref{74}, \eqref{77} and \eqref{93}, we infer that for any $t\in(0,\,+\infty)$, there exists a constant $C$ independent of $t$ such that 
	\begin{align}\label{94}
		\nonumber\int_0^t\int_Ve^{u\left(\cdot,\, t\right)}\left|\dfrac{\p u}{\p t}\right|^2 d\mu dt&=J_{\rho,h}\left(u_0\right)-J_{\rho,h}\left(u\left(\cdot,\, t\right)\right)\\
		\nonumber&\leq J_{\rho,h}\left(u_0\right)+\rho\log\left|\int_V he^{u\left(\cdot,\, t\right)}d\mu\right|-\rho\log\int_Ve^{u_0}d\mu-C\left(G,\, \rho\right)\\
		&\leq \rho\log C\left(G,\,\rho,\, h,\, u_0\right)+C\left(G,\,\rho,\, h,\, u_0\right):=C\left(G,\,\rho,\, h,\, u_0\right).
	\end{align}
	Then letting $t\rightarrow+\infty$ in \eqref{94}, one has
	$$\int_0^{+\infty}\int_Ve^{u\left(\cdot,\, t\right)}\left|\dfrac{\p u}{\p t}\right|^2 d\mu dt\leq C\left(G,\,\rho,\, h,\, u_0\right)<+\infty.$$
	Hence there exists an increasing sequence $t_n\rightarrow+\infty$ such that
	\begin{equation}\label{95}
		\lim_{n\rightarrow+\infty}\int_Ve^{u_n}\left|\dfrac{\p u_n}{\p t}\right|^2 d\mu=0.
	\end{equation}
	Here we denote $u_n=u\left(\cdot,\, t_n\right)$. In addition, $\{u_n\}$ satisfies
	\begin{equation}\label{96}
		(-\Delta)^s u_n=\rho\left(\dfrac{he^{u_n}}{\int_V he^{u_n}d\mu}-\dfrac{1}{|V|}\right)-\dfrac{\p u_n}{\p t}e^{u_n},\quad\forall\ n=1,\,2,\, \cdots.
	\end{equation}
	By means of (\ref{69}) and (\ref{95}), applying the H\"{o}lder inequality, we obtain that
	$$\int_V \left|\dfrac{\p u_n}{\p t}\right|e^{u_n}d\mu\leq\le(\int_Ve^{u_n}\left|\dfrac{\p u_n}{\p t}\right|^2 d\mu\ri)^\frac{1}{2}\le(\int_Ve^{u_0}d\mu\ri)^\frac{1}{2}\rightarrow0,\quad\textrm{as}\quad n\rightarrow+\infty,$$
	which implies that for all $x\in V$,
	\begin{equation}\label{97}
		\dfrac{\p u_n}{\p t}e^{u_n}\rightarrow0,\quad\textrm{as}\quad n\rightarrow+\infty.
	\end{equation}
	
	Now we claim that for all $n=1,\, 2,\, \cdots$, there exists a constant $C$ not depending on $t$ such that $\|u_n\|_{W^{s,2}\left(V\right)}\leq C$. Noting that a elementary inequality $e^s\geq as$ for some constant $a>0$ and all $s\in\mathbb{R}$, one has for every $x\in V$ and $n=1,\, 2,\, \cdots$,
	\begin{equation}\label{98}
		u_n\leq\frac{1}{a}e^{u_n}\leq\frac{1}{a\mu_0}\int_{V}e^{u_n}d\mu=\frac{1}{a\mu_0}\int_{V}e^{u_0}d\mu=\frac{c_0}{a\mu_0},
	\end{equation}
	where we have used \eqref{69} and the notation \eqref{146}. Thus, it follows from \eqref{98} that
	$$\int_{V}e^{2u_n}d\mu\leq e^{\frac{2c_0}{a\mu_0}}|V|:= C\left(G,\, u_0\right).$$
	Then similarly as the proof of (\ref{99}), we obtain that
	$$\left|\int_{A_{t_n}}u_nd\mu\right|\leq C\left(G,\ u_0\right),$$
	where $A_{t_n}=\big{\{}x\in V\ |\ e^{u_n}\geq2c_0/|V|\big{\}}$ for any $n=1,\, 2,\, \cdots\,$. Repeating the process of Step $3$ and Step $4$ in the proof of Lemma \ref{l10},  we conclude 
	$$\|u_n\|_{W^{s,2}\left(V\right)}\leq C\left(G,\, \rho,\, h,\, u_0\right),$$ 
	which confirms our claim. 
	
	Since $G$ is a finite garph, and all norms of the function space $C(V)$ are equivalent, we have $\|u_n\|_{\infty}\leq C$ for any $n=1,\, 2,\, \cdots$. As a consequence, passing to a subsequence $\{u_n\}$, we can look for some function $u_\infty:V\rightarrow\mathbb{R}$ such that $u_n$ converges to $u_\infty$ uniformly in $x\in V$ as $n\rightarrow+\infty$. This together with (\ref{96}) and (\ref{97}) yields that
	$$(-\Delta)^s u_\infty=\rho\left(\dfrac{he^{u_\infty}}{\int_V he^{u_\infty}d\mu}-\dfrac{1}{|V|}\right),\quad\forall\ x\in V.$$
	Thus $u_\infty$ is a solution of the mean field equation (\ref{15}).
	
	Finally, we only left to prove that along the first flow \eqref{61}, $u\left(x,\, t\right)$ converges to $u_\infty\left(x\right)$ as $t\rightarrow+\infty$ uniformly on $V$. To this end, we need an important estimate due to Lojasiewicz (\cite{Lo}, Theorem 4) in 1963. The followig finite dimensional inequality on graphs is first given by Lin-Yang (\cite{ly0}, Proposition 10), and we omit its proof.
	\begin{lemma}\label{l11}
		(Lojasiewicz-Simon inequality \cite{ly0}) Let $\sigma>0$ and $0<\theta<1/2$ be given as constants, $F\left(u\right)$ be defined by (\ref{66}), and $n$ be the number of vertexes of $V$. Along the flow (\ref{61}), if $\|u\left(\cdot,\, t\right)-u_\infty\|_{\infty}<\sigma/\sqrt{n}$ for some fixed $t$, then there exists some constant $C$ independent of $t$ such that
		\begin{equation}\label{100}
			|J_{\rho,h}\left(u\left(\cdot,\, t\right)\right)-J_{\rho,h}\left(u_\infty\right)|^{1-\theta}\leq C\left(G\right)\|F\left(u\right)\left(\cdot,\, t\right)\|_{2}.
		\end{equation}
	\end{lemma}
	Utilizing Lemma \ref{l11}, now we show that $u\left(\cdot,\, t\right)$ uniformly converges to $u_\infty$ as $t\rightarrow+\infty$.\\
	
	\noindent$\textbf{\emph{Proof of (ii) in Theorem \ref{t3}}.}$ First of all, we prove that along the flow, it always holds
	\begin{equation}\label{103}
		\lim_{t\rightarrow+\infty}\int_{V}|u\left(\cdot,\, t\right)-u_\infty|^2d\mu=0.
	\end{equation}
	Let $\sigma>0$ and $0<\theta<1/2$ be some given parameters. Observe that $u\left(x,\, t_n\right)$ converges to $u_\infty\left(x\right)$ uniformly in $x\in V$ as $n\rightarrow+\infty$. Hence for any $0<\epsilon<<\sigma$, there exists some positive integer $N>0$ such that for any $n\geq N$,
	\begin{equation}\label{101}
		\|u_n-u_\infty\|_{2}<\frac{\epsilon}{2}
	\end{equation}
	and
	\begin{equation}\label{102}
		\left(J_{\rho,h}\left(u_n\right)-J_{\rho,h}\left(u_\infty\right)\right)^\theta<\frac{\epsilon}{2}.
	\end{equation}
	Assume that
	\begin{equation}\label{104}
		t^\ast=\sup\left\{t\geq t_N\ |\ \|u\left(\cdot,\, s\right)-u_\infty\|_{2}<\sigma,\ \forall\ s\in[t_N,t]\right\}<+\infty.
	\end{equation}
	We claim that $J_{\rho,h}\left(u\left(\cdot,\, t\right)\right)>J_{\rho,h}\left(u_\infty\right)$ for all $t\in[0,\, +\infty)$. In fact, with the help of the fact that  $J_{\rho,h}\left(u\left(\cdot,\, t\right)\right)$ decreases along the flow, we can derive that $J_{\rho,h}\left(u\left(\cdot,\, t\right)\right)\geq J_{\rho,h}\left(u_\infty\right)$ for any $t\in[0,\, +\infty)$ according to \eqref{70}. If there exists some real number $t_1>0$ such that $J_{\rho,h}\left(u\left(\cdot,\ t_1\, \right)\right)=J_{\rho,h}(u_\infty)$, then we have $u\left(\cdot,\, t\right)=u_\infty$ for all $t\geq t_1$, which tells us  $t^\ast=+\infty$. However, this contradicts with (\ref{104}), showing that our claim holds. Taking into account $\p u/\p t=e^{-u\left(\cdot,\, t\right)}F\left(u\right)\left(\cdot,\, t\right)$, by means of (\ref{71}) and (\ref{100}), we infer that for any $t\in[t_N,t^\ast]$, there exsits a constant $C\left(G,\, \theta\right)>0$ such that
	\begin{align}\label{105}
		\nonumber	-\dfrac{d}{dt}\left(J_{\rho,h}\left(u\left(\cdot,\, t\right)\right)-J_{\rho,h}\left(u_\infty\right)\right)^\theta&=\theta\left(J_{\rho,h}\left(u\left(\cdot,\, t\right)\right)-J_{\rho,h}\left(u_\infty\right)\right)^{\theta-1}\int_VF\left(u\right)\left(\cdot,\, t\right)\dfrac{\p u}{\p t} d\mu\\
		\nonumber&\geq C\left(G,\, \theta\right)\frac{\int_Ve^{-u\left(\cdot,\, t\right)}F^2(u)\left(\cdot,\, t\right) d\mu}{\|F\left(u\right)\left(\cdot,\, t\right)\|_{2}}\\
		\nonumber&\geq C\left(G,\, \theta\right)e^{-C}\|F(u(\cdot,\,t))\|_2\\
		&\geq C\left(G,\, \theta\right)\left\|\dfrac{\p u}{\p t}\right\|_{2},
	\end{align}
	where we have used the fact that $\|u\left(\cdot,\, t\right)\|_{\infty}\leq C$ for any $t\in[t_N,t^\ast]$, in addition, our assertion $J_{\rho,h}\left(u\left(\cdot,\, t\right)\right)>J_{\rho,h}\left(u_\infty\right)$ ensures that the first equality of \eqref{105} holds for any $t\in[t_N,t^\ast]$. Integrating the both sides of (\ref{105}) from $t_N$ to $t^\ast$, from (\ref{102}), we derive that
	\begin{equation}\label{106}
		\int_{t_N}^{t^\ast}\left\|\dfrac{\p u}{\p t}\right\|_{2}dt\leq C\left(J_{\rho,h}\left(u_N\right)-J_{\rho,h}\left(u_\infty\right)\right)^\theta<\frac{\epsilon}{2}.
	\end{equation}		
	Moreover, it follows from the H\"{o}lder inequality that for any $t\in[t_N,t^\ast]$,
	\begin{equation}\label{107}
		\dfrac{d}{dt}\left(\int_V |u\left(\cdot,\, t\right)-u_\infty|^2d\mu\right)^{\frac{1}{2}}=\dfrac{1}{\|u\left(\cdot,\, t\right)-u_\infty\|_{2}}\int_V\left|u\left(\cdot,\, t\right)-u_\infty\right|\bigg{|}\dfrac{\p u}{\p t}\bigg{|}d\mu\leq\left\|\dfrac{\p u}{\p t}\right\|_{2}.
	\end{equation}
	Integrating (\ref{107}) from $t_N$ to $t^\ast$, we conclude from (\ref{101}) and (\ref{106}) that
	\begin{equation*}
		\|u\left(\cdot,\, t^\ast\right)-u_\infty\|_{2}\leq\int_{t_N}^{t^\ast}\left\|\dfrac{\p u}{\p t}\right\|_{2}dt+\|u_N-u_\infty\|_{2}<\epsilon<<\sigma,
	\end{equation*}
	in contradiction with the definition of $t^\ast$. Thus we have  $t^\ast=+\infty$ and $\|u\left(\cdot,\, t\right)-u_\infty\|_{2}<\sigma$ for all $t\geq t_N$, i.e, \eqref{103} holds.
	
	Finally, it follows from \eqref{103} that
	$$\lim_{t\rightarrow+\infty}\sum_{i=1}^n\mu\left(x_i\right)|u\left(x_i,\, t\right)-u_\infty\left(x_i\right)|^2=0.$$
	For every $i=1,\, 2,\, \cdots,\ n$, we have
	\begin{equation*}
		|u\left(x_i,\, t\right)-u_\infty\left(x_i\right)|^2\leq\frac{1}{\mu_0}\sum_{i=1}^n\mu\left(x_i\right)|u\left(x_i,\, t\right)-u_\infty\left(x_i\right)|^2\rightarrow0,\quad\textrm{as}\quad t\rightarrow+\infty.
	\end{equation*}
	As a consequence, along the flow (\ref{61}), $u\left(x,\, t\right)\rightarrow u_\infty\left(x\right)$ uniformly on $V$ as $t\rightarrow+\infty$. This completely ends the proof of Theorem \ref{t3}.$\hfill\Box$
	
	\subsection{The second flow}
	
	Recall the second heat flow on a finite graph $G=\left(V,\ E\right)$, which is defined by
	\begin{equation}\label{63}\left\{\begin{array}{lll}
			\dfrac{\p}{\p t}u=\alpha(t)h-R,\quad\left(x,\, t\right)\,\in \, V\times\left(0,\, +\infty\right),\\[3ex]
			u\left(x,\, 0\right)=u_0\left(x\right)\in\mathscr{M},\quad\ \  x\,\in \,V.\end{array}\right.
	\end{equation}
	Here $\rho$ is a positive constant, $\mathscr{M}$ is given by \eqref{35}, the scalar curvature type function $R$ and the normalized coefficient $\alpha(t)$ are defined by
	\begin{equation}\label{64}
		R=\frac{1}{\rho}e^{-u}\le((-\D)^su+\frac{\rho}{|V|}\ri),\quad	\alpha(t)=\frac{\int_{V}hRe^ud\mu}{\int_{V}h^2e^ud\mu}.
	\end{equation}
	
	Due to the initial function $u_0\in\mathscr{M}$, we first have the following obsevation. 
	\begin{lemma}\label{l12}
		For all $t$ whenever the flow \eqref{63} exists, there always holds
		\begin{equation}\label{110}
			\int_Vhe^{u(\cdot,\,t)}d\mu=\int_Vhe^{u_0}d\mu=1.
		\end{equation}
	\end{lemma}
	\begin{proof}
		Differentiating $\int_Vhe^{u(\cdot,\,t)}d\mu$ with respect to $t$, in view of \eqref{63} and \eqref{64}, one has
		\begin{align*}
			\dfrac{d}{d t}\int_Vhe^{u(\cdot,\,t)}d\mu&=\int_Vhe^{u(\cdot,\,t)}\dfrac{\p u}{\p t}d\mu\\
			&=\int_Vhe^{u(\cdot,\,t)}(\alpha(t)h-R)d\mu\\
			&=\alpha(t)\int_Vh^2e^{u(\cdot,\,t)}d\mu-\int_VhRe^{u(\cdot,\,t)}d\mu\\
			&=0,
		\end{align*}
		which together with $u_0\in\mathscr{M}$ yields that \eqref{110} holds.
	\end{proof}
	By Lemma \ref{l12}, the energy functional corresponding to \eqref{15} reduces to
	\begin{equation}\label{111}
		J\left(u\right)=\frac{1}{2}\int_V |\nabla^s u|^2d\mu+\frac{\rho}{|V|}\int_V ud\mu,\quad\forall\ u\in\mathscr{M}.
	\end{equation}
	Next we shall verify that $J(u)$ is non-increasing along the second flow \eqref{63} for any $t\geq0$ if the flow exists.
	\begin{lemma}\label{l13}
		Suppose the flow exits on $[0,\,T)$ ($T$ may be $+\infty$). For any $0\leq t_1<t_2<T$, we have
		\begin{equation}\label{112}
			J\left(u\left(\cdot,\,t_2\right)\right)\leq J\left(u\left(\cdot,\,t_1\right)\right).
		\end{equation}
	\end{lemma}
	\begin{proof}
		Differentiating \eqref{111} with respect to $t$, by Lemma \ref{l2}, we have from \eqref{63}, \eqref{64}
		\begin{align}\label{114}
			\nonumber\dfrac{d}{d t}J(u(\cdot,\,t))&=\int_V\le((-\D)^su+\frac{\rho}{|V|}\ri)\dfrac{\p u}{\p t}d\mu\\
			\nonumber&=-\rho\int_VR^2e^ud\mu+\rho\alpha(t)\int_VhRe^ud\mu\\
			\nonumber&=-\rho\int_VR^2e^ud\mu+2\rho\alpha(t)\int_VhRe^ud\mu-\rho\alpha^2(t)\int_{V}h^2e^ud\mu\\
			&=-\rho\int_V\le(\alpha(t)h-R\ri)^2e^ud\mu\leq0.
		\end{align}
		Hence the lemma follows from this immediately.
	\end{proof}
	
	\subsubsection{Short and long time existence}   
	The short time existence follows from the standard ordinary differential equation theory, which we omit its details for
	saving the length of the paper. Let us define
	\begin{equation}\label{116}
		T=\sup\Big{\{}T_0>0\ |\ u:V\times[0,\, T_0]\rightarrow\mathbb{R}\ \mathrm{is\ the\ unique\ solution\ of\ the\ second\ flow}\Big{\}}<+\infty.
	\end{equation}
	In what follows, we mainly discuss the long time existence of the second flow. Making use of Lemma \ref{l13}, it follows from \eqref{114} that for a finite $T>0$
	\begin{equation}\label{115}
		J(u(T))+\rho\int_0^T\int_V\le(\alpha(t)h-R\ri)^2e^ud\mu dt=J(u_0).
	\end{equation}
	Henceforth we often denote $u(\cdot,\,t)$ by $u(t)$ for simplicity. Applying Lemma $\ref{l4}$, for any $\beta>0$ and $t\in[0,\, T)$,  there exists a constant $C$ such that
	\begin{equation}\label{117}
		\int_V e^{\beta\frac{\left(u\left(t\right)-\overline{u}\left(t\right)\right)^2}{\|\nabla^s u\left(t\right)\|^2_{2}}}d\mu\leq C\left(G,\, \beta\right).
	\end{equation}
	Here $\bar{u}\left(t\right)$ denotes the integral average of $u(t)$ over $V$. By using (\ref{117}) and the Young inequality, noting that $\rho>0$, we have for any $\epsilon>0$ and $t\in[0,\, T)$
	\begin{equation}\label{118}
		\frac{\rho}{|V|}\int_{V}u\left(t\right)d\mu\geq\rho\log\int_Ve^{u\left(t\right)}d\mu-\epsilon\rho\int_{V}|\nabla^s u\left(t\right)|^2d\mu+C\left(G,\, \epsilon,\, \rho\right).
	\end{equation}
	Let us insert (\ref{118}) into (\ref{111}) and take $\epsilon=1/\left(2\rho\right)$. In view of (\ref{112}), one can obtain that for any $t\in[0,\, T)$,
	\begin{align}\label{119}
		J\left(u_0\right) \geq J\left(u\left(t\right)\right)
		\geq\rho\log\int_Ve^{u\left(t\right)}d\mu+C\left(G,\, \rho\right),
	\end{align}
	which leads to
	\begin{equation}\label{121}
		\int_Ve^{u\left(t\right)}d\mu\leq C(G,\,\rho)e^{\frac{1}{\rho}J(u_0)}:=C(G,\,\rho,\,u_0).
	\end{equation}
	On the other hand, it follows from \eqref{110} that
	$$1=\int_Vhe^{u(t)}d\mu\leq\max_{V}h\int_Ve^{u\left(t\right)}d\mu.$$
	Consequently, it results for any $t\in[0,\, T)$
	\begin{equation}\label{120}
		\int_Ve^{u\left(t\right)}d\mu\geq\frac{1}{\max_{V}h}>0.
	\end{equation}
	Combining \eqref{119} and \eqref{120}, we deduce that for any $t\in[0,\, T)$
	$$J\left(u\left(t\right)\right)\geq-\rho\log(\max_{V}h)+C\left(G,\, \rho\right).$$
	If the flow \eqref{63} globally exists in a time interval $[0,\,+\infty)$, then by letting $T\rightarrow+\infty$ in \eqref{115}, one arrives at
	\begin{equation}\label{125}
		\int_0^{+\infty}\int_V\le(\alpha(t)h-R\ri)^2e^ud\mu dt\leq\log(\max_{V}h)+C\left(G,\, \rho\right)+\frac{1}{\rho}J(u_0):=C(G,\,\rho,\,h,\,u_0).
	\end{equation}
	
	In order to get the upper bound for $\|u(t)\|_{W^{s,2}(V)}$, we have the following estimate.
	\begin{lemma}\label{l14}
		For any $\gamma\geq1$, there exists a constant $C$ depending only on the graph $G$, $T$, $\rho$, $h$, $u_0$ and $\gamma$, such that the sloution $u(t)$ of the flow \eqref{63} which exists on $[0,\,T)$ satisfies
		\begin{equation*}
			\sup_{0\leq t<T}\int_Ve^{\gamma u(t)}d\mu\leq C(G,\,T,\,\rho,\,h,\,u_0,\,\gamma).
		\end{equation*}
	\end{lemma}
	\begin{proof}
		With the help of \eqref{64} and \eqref{110}, we obtain for any $t\in[0,\, T)$
		\begin{align}\label{122}
			\nonumber\int_V(\alpha(t)h-R)e^{u(t)}d\mu&=\alpha(t)\int_Vhe^{u(t)}d\mu-\int_VRe^{u(t)}d\mu\\
			\nonumber&=\alpha(t)-\int_V\frac{1}{\rho}\le((-\D)^su+\frac{\rho}{|V|}\ri)d\mu\\
			&=\alpha(t)-1.
		\end{align}
		Making use of the H\"{o}lder inequality, by \eqref{121} and \eqref{122}, we derive
		\begin{align}\label{123}
			\nonumber(\alpha(t)-1)^2&=\le(\int_V(\alpha(t)h-R)e^{\frac{u(t)}{2}}e^{\frac{u(t)}{2}}d\mu\ri)^2\\
			\nonumber&\leq\int_V(\alpha(t)h-R)^2e^{u(t)}d\mu\int_Ve^{u\left(t\right)}d\mu\\
			&\leq C(G,\,\rho,\,u_0)\int_V(\alpha(t)h-R)^2e^{u(t)}d\mu.
		\end{align} 
		Observe the elementary inequality $|\alpha(t)|\leq|\alpha(t)-1|+1\leq2((\alpha(t)-1)^2+1)$. Then an easy calculation yields that
		\begin{align}\label{124}
			\nonumber\dfrac{d}{d t}\int_Ve^{\gamma u(t)}d\mu&=\gamma\int_Ve^{\gamma u(t)}(\alpha(t)h-R)d\mu\\
			\nonumber&=-\frac{\gamma}{\rho}\int_Ve^{(\gamma-1) u(t)}(-\D)^su(t)d\mu-\frac{\gamma}{|V|}\int_Ve^{(\gamma-1) u(t)}d\mu+\alpha(t)\gamma\int_Ve^{\gamma u(t)}hd\mu\\
			&\leq 2\gamma\max_{V}h((\alpha(t)-1)^2+1)\int_Ve^{\gamma u(t)}d\mu,
		\end{align}
		where we use the fact that
		\begin{align*}
			-\frac{\gamma}{\rho}\int_Ve^{(\gamma-1) u(t)}(-\D)^su(t)d\mu&=-\frac{\gamma}{\rho}\int_V\nabla^se^{(\gamma-1) u(t)}\nabla^su(t)d\mu\\
			&=-\frac{\gamma}{2\rho}\sum_{x\in V}\sum_{y\in V, y\not=x}W_s(x,y)\left(e^{(\gamma-1)u\left(x,\, t\right)}-e^{(\gamma-1)u\left(y,\, t\right)}\right)\left(u\left(x,\,t\right)-u\left(y,\, t\right)\right)\\
			&\leq0.
		\end{align*}
		By means of \eqref{125} and \eqref{123}, it follows from \eqref{124} that for any $t\in[0,\, T)$
		\begin{align*}
			\log\int_Ve^{\gamma u(t)}d\mu-\log\int_Ve^{\gamma u_0}d\mu&\leq C(\gamma,\,h)\le(C(G,\,\rho,\,u_0)\int_0^t\int_V(\alpha(t)h-R)^2e^{u(t)}d\mu dt+t\ri)\\
			&\leq C(\gamma,\,h)\le(C(G,\,\rho,\,h,\,u_0)+T\ri),
		\end{align*}
		which implies that
		$$\int_Ve^{\gamma u(t)}d\mu\leq e^{C(\gamma,\,h)\le(C(G,\,\rho,\,h,\,u_0)+T\ri)}\int_Ve^{\gamma u_0}d\mu:=C(G,\,T,\,\rho,\,h,\,u_0,\,\gamma).$$
		Hence we conclude
		\begin{equation}\label{113}
			\sup_{0\leq t<T}\int_Ve^{\gamma u(t)}d\mu\leq C(G,\,T,\,\rho,\,h,\,u_0,\,\gamma).
		\end{equation}
		This ends the proof of the lemma.
	\end{proof}
	In view of Lemma \ref{l14}, now we are able to state the following norm estimate of the solution $u(t)$ for the flow \eqref{63}. 
	\begin{lemma}\label{l15}
		For any $t\in[0,\,T)$, if $u\left(t\right)$ is the solution of the second flow \eqref{63}, then there exists a constant $C$ depending only on $G$, $T$, $\rho$, $h$ and $u_0$ such that
		\begin{equation*}
			\|u\left(t\right)\|_{W^{s,2}\left(V\right)}\leq C\left(G,\, T,\, \rho,\, h,\, u_0\right),\quad\forall\, t\in[0,\,T).
		\end{equation*}
	\end{lemma}
	\begin{proof}
		From \eqref{121} and \eqref{120}, we obtain that for any $t\in[0,\, T)$
		\begin{equation}\label{127}
			0<\frac{1}{\max_{V}h}\leq\int_Ve^{u\left(t\right)}d\mu\leq C(G,\,\rho,\,u_0).
		\end{equation}
		We define a set
		\begin{equation}\label{128}
			A_{t}=\left\{x\in V\ |\ e^{u\left(t\right)}\geq\frac{1}{2|V|\max_{V}h}\right\}\not=\varnothing.
		\end{equation}
		It follows from (\ref{113}), \eqref{127} and (\ref{128}) that
		\begin{align*}
			\nonumber		\frac{1}{\max_{V}h}&\leq\int_Ve^{u\left(t\right)}d\mu\\
			\nonumber	 &=\int_{V\setminus A_{t}} e^{u\left(t\right)}d\mu+\int_{A_{t}} e^{u\left(t\right)}d\mu\\
			\nonumber			&\leq\frac{1}{2|V|\max_{V}h}\int_{V\setminus A_{t}}1 d\mu+\left(\int_{V} e^{2u\left(t\right)}d\mu\right)^{\frac{1}{2}}\left(\int_{A_t} 1d\mu\right)^{\frac{1}{2}}\\
			&\leq\frac{1}{2\max_{V}h}+C\left(G,\, T,\, \rho,\, h,\, u_0\right)|A_t|^{\frac{1}{2}}.
		\end{align*}
		Consequently, it yields that
		\begin{equation}\label{130}
			|A_t|\geq C^{-2}\left(G,\, T,\, \rho,\, h,\, u_0\right)\frac{1}{4(\max_{V}h)^2}:=C\left(G,\, T,\, \rho,\, h,\, u_0\right)>0.
		\end{equation}
		Utilizing Jensen's inequality, we infer from \eqref{128} that
		\begin{equation}\label{131}
			\le(-\log(2|V|\max_{V}h)\ri)|A_t|\leq\int_{A_t} u\left(t\right)d\mu\leq |A_t|\log\le(|A_t|^{-1}\int_{A_t}e^{u\left(t\right)}d\mu\ri).
		\end{equation}
		Combining \eqref{127}, \eqref{130} and \eqref{131}, one has
		\begin{equation}\label{132}
			\bigg{|}\int_{A_t} u\left(t\right)d\mu\bigg{|}\leq C\left(G,\, T,\, \rho,\, h,\,u_0\right).
		\end{equation}
		Hence, for any $\epsilon>0$, it follows from \eqref{132} that
		\begin{align}\label{133}
			\nonumber	\left(\int_V u\left(t\right)d\mu \right)^2
			\leq&\left(\int_{A_t} u\left(t\right)d\mu\right)^2+\left(\int_{V\setminus A_t} u\left(t\right)d\mu\right)^2+\frac{1}{\epsilon}\left(\int_{A_t} u\left(t\right)d\mu\right)^2
			+\epsilon\left(\int_{V\setminus A_t} u\left(t\right)d\mu\right)^2\\
			\leq&\left(1+\frac{1}{\epsilon}\right)C\left(G,\, T,\, \rho,\, h,\, u_0\right)+\left(1+\epsilon\right)|V\setminus A_t|\int_V u^2\left(t\right)d\mu.
		\end{align}	
		By Lemma \ref{l3}, in view of \eqref{133}, we deduce that
		\begin{align}\label{134}
			\nonumber\int_V u^2\left(t\right)d\mu\leq& C(G)\int_V |\nabla^s u\left(t\right)|^2d\mu+\int_V \overline{u}^2\left(t\right)d\mu\\
			\leq& C(G)\int_V |\nabla^s u\left(t\right)|^2d\mu+\dfrac{\left(1+\epsilon\right)|V\setminus A_t|}{|V|}\int_V u^2\left(t\right)d\mu+C(G,\, T,\, \rho,\, h,\, u_0,\,\epsilon).
		\end{align}
		Choose small $\epsilon>0$ in \eqref{134} so that
		\begin{equation}\label{135}
			\int_V u^2\left(t\right)d\mu\leq C\left(G\right)\int_{V}|\nabla^s u\left(t\right)|^2d\mu+C\left(G,\, T,\, \rho,\, h,\, u_0\right),\quad\forall\, t\in[0,\, T).
		\end{equation}
		According to \eqref{112}, applying the Young inequality to \eqref{111} implies that for any $t\in[0,\, T)$ and $\epsilon>0$,
		\begin{equation}\label{136}
			J\left(u_0\right)\geq J\left(u(t)\right)\geq\frac{1}{2}\int_V |\nabla^s u\left(t\right)|^2d\mu-\frac{\rho\epsilon}{|V|}\int_V u^2\left(t\right)d\mu-\frac{\rho}{4\epsilon}.
		\end{equation}
		Inserting (\ref{135}) into (\ref{136}), we obtain
		$$\int_V |\nabla^s u\left(t\right)|^2d\mu-\frac{\rho}{2\epsilon}-2J(u_0)\leq\frac{2\rho\epsilon}{|V|}\int_V u^2\left(t\right)d\mu\leq\frac{2\rho \epsilon}{|V|}\left(C(G)\int_{V}|\nabla^s u\left(t\right)|^2d\mu+C\left(G,\, T,\, \rho,\, h,\, u_0\right)\right).$$
		Let us take $\epsilon=|V|/\left(4\rho C(G)\right)$, we have
		\begin{equation}\label{137}
			\int_V |\nabla^s u\left(t\right)|^2d\mu\leq C\left(G,\, T,\, \rho,\, h,\, u_0\right).
		\end{equation}
		Noting that \eqref{135} and (\ref{137}), we get
		$$\|u\left(t\right)\|_{W^{s,2}\left(V\right)}\leq C\left(G,\, T,\, \rho,\, h,\, u_0\right),\quad\forall\ t\in[0,\,T).$$
		This completes the proof.
	\end{proof}
	
	Now we would like to give the long time existence of the flow \eqref{63}.\\
	
	\noindent$\textbf{\emph{Proof of (i) in Theorem \ref{t4}}.}$ Let $T$ be given in (\ref{116}) satisfying $T<+\infty$. By Lemma \ref{l15}, it is not hard to see that $u\left(t\right)$ can be extended beyond $T$, contradicting to (\ref{116}). Thus $T=+\infty$. Hence the solution $u(t)$ of the second flow \eqref{63} with given initial function $u_0\in\mathscr{M}$ must exist globally. $\hfill\Box$
	
	\subsubsection{Global convergence}
	
	Finally, we are devoted to obtaining the convergence results of the solution $u(t)$ for the flow \eqref{63} in this part. At first it follows from \eqref{63} and \eqref{125} that
	$$\int_0^{+\infty}\int_Ve^{u(t)}\le|\dfrac{\p u}{\p t}\ri|^2d\mu dt\leq C(G,\,\rho,\,h,\,u_0)<+\infty.$$
	Therefore there is a positive sequence $\{t_n\}$ satisfying $t_n<t_{n+1}$ and $t_n\rightarrow+\infty$ as $n\rightarrow+\infty$ such that
	\begin{equation}\label{138}
		\lim_{n\rightarrow+\infty}\int_Ve^{u\left(t_n\right)}\left|\dfrac{\p u\left(t_n\right)}{\p t}\right|^2 d\mu=0.
	\end{equation}
	Let us write $u_n=u\left(t_n\right)$. In view of (\ref{121}) and (\ref{138}), we obtain that
	\begin{align*}
		\int_V \left|\dfrac{\p u_n}{\p t}\right|e^{u_n}d\mu&\leq\le(\int_Ve^{u_n}\left|\dfrac{\p u_n}{\p t}\right|^2 d\mu\ri)^\frac{1}{2}\le(\int_Ve^{u_n}d\mu\ri)^\frac{1}{2}\\
		&\leq C(G,\,\rho,\,u_0)\le(\int_Ve^{u_n}\left|\dfrac{\p u_n}{\p t}\right|^2 d\mu\ri)^\frac{1}{2}\rightarrow0,\quad\textrm{as}\quad n\rightarrow+\infty.
	\end{align*}
	Hnece we have for every $x\in V$,
	\begin{equation}\label{139}
		\dfrac{\p u_n}{\p t}e^{u_n}\rightarrow0,\quad\textrm{as}\quad n\rightarrow+\infty.
	\end{equation}
	
	Next similarly as the proof of \eqref{98}, we obtain from \eqref{121} that for all $x\in V$ and $n=1,\, 2,\, \cdots$,
	\begin{equation}\label{140}
		u_n\leq e^{u_n}\leq\frac{1}{\mu_0}\int_{V}e^{u_n}d\mu\leq C\left(G,\,\rho,\, u_0\right),
	\end{equation}
	Then for any $\gamma\geq1$, it follows from \eqref{140} that for all $n=1,\, 2,\, \cdots$,
	$$\int_{V}e^{\gamma u_n}d\mu\leq e^{\gamma C\left(G,\,\rho,\, u_0\right)}|V|:= C\left(G,\,\rho,\, u_0,\,\gamma\right).$$
	Repeating the process of Lemma \ref{l15}, we conclude that  there exists a constant $C$ not depending on $t$ such that for any $n=1,\, 2,\, \cdots$, 
	\begin{equation}\label{141}
		\|u_n\|_{W^{s,2}\left(V\right)}\leq C\left(G,\, \rho,\, h,\, u_0\right).
	\end{equation}
	
	Finally, we shall show the global convergence of $u_n$ as $n\rightarrow+\infty$. \\
	
	\noindent$\textbf{\emph{Proof of (ii) in Theorem \ref{t4}}.}$ In view of \eqref{63}, \eqref{64} and \eqref{110}, we derive that for any $x\in V$ and $n=1,\, 2,\, \cdots$, $\{u_n\}$ satisfies
	\begin{equation}\label{142}
		\rho e^{u_n}\dfrac{\p u_n}{\p t}=\alpha(t_n)\rho h e^{u_n}-\le((-\D)^su_n+\frac{\rho}{|V|}\ri),
	\end{equation}
	provided with $u_n\in\mathscr{M}$. In addition, with \eqref{63} and \eqref{123} at hand, one has from \eqref{138} that
	\begin{align*}
		0\leq(\alpha(t_n)-1)^2&\leq C(G,\,\rho,\,u_0)\int_V(\alpha(t_n)h-R)^2e^{u_n}d\mu\\
		&=C(G,\,\rho,\,u_0)\int_V\le|\dfrac{\p u_n}{\p t}\ri|^2e^{u_n}d\mu\rightarrow 0,\quad\textrm{as}\quad n\rightarrow+\infty.
	\end{align*}
	Thus we infer that
	\begin{equation}\label{126}
		\lim_{n\rightarrow+\infty}\alpha(t_n)=1.
	\end{equation}
	Due to the equivalent norms, it follows from \eqref{141} that  $\|u_n\|_{\infty}\leq C$ for every $n=1,\, 2,\, \cdots$. Hence, there exists a subsequence of $\{u_n\}$ and some function $u_\infty$ such that $u_n\rightarrow u_\infty$ uniformly in $x\in V$ as $n\rightarrow+\infty$. Moreover, it is possible to see that $u_\infty\in\mathscr{M}$. Let $n$ tends to $+\infty$ in \eqref{142}. By means of \eqref{139} and \eqref{126}, we deduce that $u_\infty$ satisfies
	$$(-\Delta)^s u_\infty=\rho he^{u_\infty}-\dfrac{\rho}{|V|},\quad\forall\ x\in V.$$
	Consequently, $u_\infty$ is a solution of the mean field equation (\ref{15}), which completely ends the proof of our last theorem.   $\hfill\Box$
	
	\begin{remark}
		(i) Similar to the fist flow \eqref{61}, we may also prove that the solution $u(t)$ of flow \eqref{63} uniformly converges to $u_\infty$ on $V$ when $t$ tends to $+\infty$ (details are left
		for the reader).\\
		(ii) As we might see, by using different methods, we have obtained the solutions of \eqref{15} in three cases. However, it is currently unclear whether these solutions are the same or different. 
 	\end{remark}
	
	\noindent\textbf{Conflict of interest $\&$ Data availability}
	The author declares that there are no conflicts of interests regarding the publication of this paper. Data sharing not applicable as no datasets were used or analysed during the current study.

\end{document}